\documentclass[11pt]{amsart}
\numberwithin{equation}{section} \oddsidemargin=-.0cm
\usepackage{amsfonts, amssymb, amsmath, amsthm}
\usepackage[colorlinks=true, pdfstartview=FitV, linkcolor=blue,
            citecolor=blue, urlcolor=blue]{hyperref}
\usepackage{pifont}

\def\M{\mathfrak{M}}
\def \no#1#2#3 {{\bf #1} (#3), #2.}
\def \eds#1#2#3 {#1, #2, #3.}

\newtheorem{proposition}{Proposition}[section]
\newtheorem{theorem}[proposition]{Theorem}

\theoremstyle{definition}
\newtheorem{definition}[proposition]{Definition}
\newtheorem{remark}[proposition]{Remark}

\newtheorem*{Aonr}{Assumptions on the function $r$}
\newtheorem*{CrR}{Cross condition}
\numberwithin{equation}{section}
\usepackage{color}

\definecolor{rred}{rgb}{0.7,0.,0.}
\newcommand{\mk}{\color{black}}
\newcommand{\mkk}{\color{black}}
\newcommand{\jr}{\color{black}}
\newcommand{\mko}{\color{black}}
\title[Conjugacy problems and the Koopman operator]{Homeomorphisms group of normed vector space:\\ Conjugacy problems and the
Koopman operator}

\author[M. D. Chekroun and J. Roux] {Micka\"el D. Chekroun and Jean Roux}
\date{}
{\mkk \address[MDC]{Department of Mathematics, University of Hawaii at Manoa, Honolulu, HI 96822, USA;  and 
Department of Atmospheric \& Oceanic Sciences, University of California, Los Angeles, CA 90095-1565, USA} \email{mchekroun@atmos.ucla.edu}}

{\mkk \address[JR]{CERES-ERTI, \'Ecole Normale
Sup\'erieure, 75005 Paris, France}
\email{jroux@lmd.ens.fr}}



\begin{document}

\maketitle
\begin{abstract}
This article is concerned with conjugacy problems arising in {\jr the} 
homeomorphisms group, Hom($F$), of {\mk unbounded} subsets $F$ of normed
vector spaces $E$. Given two homeomorphisms $f$ and $g$ in Hom($F$),
it is shown how the existence of a conjugacy may be related to the
existence of a common generalized eigenfunction of the associated
Koopman operators. This common eigenfunction serves to build a
topology on Hom($F$), where the conjugacy is obtained as limit of a
sequence generated by the conjugacy operator, when this limit
exists. The main conjugacy theorem is presented in a class of
generalized Lipeomorphisms.
\end{abstract}

\section{Introduction}
In this article we consider the conjugacy problem in the
homeomorphisms group of a finite dimensional normed vector space
$E$. It is out of the scope of the present work to review the
problem of conjugacy in general, and  the reader may consult for
instance \cite{TOPOChaperon86, TOPOCong, TOPOkat,TOPOkotus82,
TOPOIrwin, TOPOPalis, TOPORou, TOPOSmale0,TOPOSmale} and references therein, to get a
partial survey  of the question from a dynamical point of view. The
present work raises the problem of conjugacy in the group Hom($F$)
{\jr consisting of}  homeomorphisms of an {\mk unbounded} subset $F$ of $E$ and
is intended to demonstrate how the conjugacy problem, in such a
case, may be related to spectral properties of the associated
Koopman operators. In this sense, this paper provides new insights
on the relations between the {\mkk spectral theory of  dynamical
systems \cite{TOPO_V00,TOPOSinai_book,TOPO_goodson, LM94}} and the {\mkk topological conjugacy problem \cite{TOPOSmale0, TOPOSmale}}\footnote{{\mkk Usually, the spectral theory of  dynamical
systems makes usage of concepts from ergodic theory, but the latter are not required in the description of the relationships that we propose below.}}.

More specifically, given two homeomorphisms $f$ and $g$ of $F$, we
show here that the conjugacy problem in Hom($F$) is {\jr related}  to the
existence  of a common generalized eigenfunction for the associated
Koopman operators $U_f$ and $U_g$ (cf. Definition \ref{Def_Koop}),
{\it i.e.} a function $\Phi$ satisfying,
\begin{equation}\label{System_cohomo_intro} \left \{
\begin{array}{l}
U_f(\Phi)\geq  \lambda  \Phi, \\
U_g(\Phi) \geq  \mu \Phi,\\
\end{array} \right.
\end{equation}
for some $\lambda, \mu >0$, where $\Phi$ lives within some cone $K$
of the set of continuous real-valued functions on $F$. The elements
of this cone {\mk possess} the particularity {\jr of exhibiting}  a behavior at
infinity prescribed by a subadditive function $R$; see Section
\ref{Sec_functional-setting}.

More precisely, when such a $ \Phi$ exists,  it is shown how $\Phi$
can be used to build a topology such that the sequence of iterates
$\{\mathcal{L}_{f,g}^n(h_0)\}_{n\in \mathbb{N}}$, of the conjugacy
operator\footnote{where  {\mkk for $f$ and $g$ given in Hom$(F)$,} $\mathcal{L}_{f,g}:\psi \mapsto f\circ \psi
\circ g^{-1}$ is acting on $\psi \in $Hom$(F)$.} initiated to some
$h_0 \in$ Hom($F$) close enough to $\mathcal{L}_{f,g}(h_0)$ in that
topology, converges to the conjugacy $h$ satisfying $f\circ h=h\circ
g$, provided that $\{\mathcal{L}_{f,g}^n(h_0)\}$ is bounded on every
compact {\mk subset} of $F$; cf. Theorem \ref{THM_conjug}. The topology built
from $\Phi$ relies on a premetric on Hom($F$) where $\Phi$ serves
to {\jr weigh} the distance to the identity of any homeomorphism of $F$; {\mkk see Eqs. \eqref{Eq_def-semi-norm} and \eqref{Eq_def_rho} below}.

The plan of this article is as follows. Section
\ref{Sec_functional-setting} sets up the functional framework used
in this article, where in particular the main properties of the
topology built from any member $\Phi \in K$ are derived with 
particular attention to closure properties and convergence in that
topology of sequences in Hom($F$); cf Propositions
\ref{Thm_Closure1} and \ref{Thm_Closure2}.  Section
\ref{Sec_Fixed_Pt_Thm} establishes a fixed point theorem, Theorem
\ref{Thm_fixed-point}, for mappings acting on Hom($F$), when this
{\jr group} is endowed with the topology discussed in Section
\ref{Sec_functional-setting}. In section \ref{Sec_main_thm} the main
theorem of conjugacy, Theorem \ref{THM_conjug}, is proved based on
Theorem \ref{Thm_fixed-point} applied to the conjugacy operator,
where the contraction property is shown to be {\jr related to the existence of a common generalized eigenfunction $\Phi$ of  a
generalized eigenvalue problem of type \eqref{System_cohomo_intro}}.
This related generalized eigenvalue problem for the Koopman
operators associated with the conjugacy problem is then discussed in
Section \ref{SEC_COHOM_discussion} where{\jr ,} in particular{\jr ,} connections
with relatively recent results about {\mk functional equations such as the} Schr\"oder equations and {\mk the} Abel
equation are established.  {\mk Concluding remarks regarding the possible extensions of the present work are presented in Section 5.}
The results obtained {\mk in the present study} were motivated
in part by \cite{TOPOChek_al06}, where {\mkk results are derived 
for the conjugacy problem on not necessarily compact manifolds}. {\mk Connections with topological equivalence
 problems between periodic vector fields and autonomous ones as considered in \cite{TOPOChek_al06},  will be discussed elsewhere. }

\section{A functional framework on the homeomorphisms group}\label{Sec_functional-setting}
In this section we introduce {\jr a} family of subgroups of homeomorphisms
for the composition law. These subgroups associated with the
framework from which they are derived, will be used in the analysis
of the conjugacy problem in the homeomorphisms group itself. The
topology with which they are endowed is  introduced here and the
main properties are derived. The extension of these topologies to
the whole group of homeomorphisms is also presented and the related
closure properties and convergence of sequences in the
homeomorphisms group are discussed.

\subsection{Notations and preliminaries}\label{Sec_prellim1}
In this article $E$ denotes a $d$-dimensional normed vector space
($d \in \mathbb{N}^*$), endowed with a norm denoted by $\|\cdot\|$
and $F$ denotes an {\mk unbounded} subset of $E$. The following class of
functions serves to specify some behavior at infinity of
homeomorphisms and to build topologies that will be central in our
approach; cf Proposition \ref{Thm_Main-propertiesI}.

\begin{definition}\label{Def_SpaceE_F}
{\bf The space $\mathcal{E}_F^R$.} Let $R:\mathbb{R}^+ \rightarrow
\mathbb{R}^{+}-\{0\},$ be a continuous function, which is
subadditive, {\it i.e.}, $$R(u+v)\leq R(u)+R(v), \; \forall \;
u,v\in \mathbb{R}^{+}.
$$

We denote by $\mathcal{E}^R_F$ the set of continuous functions
$\Phi:F\rightarrow \mathbb{R}^+$, satisfying:
\begin{itemize}
\item[(G$_1$)] $\exists \; m>0,\; \forall \; x\in F,\; m\leq \Phi(x),$
\item[(G$_2$)] Coercivity condition: $\Phi(x) \longrightarrow +\infty$, as $x \in F$ and
$\|x\| \rightarrow +\infty$ ,
\item[(G$_3$)] Cone condition: There exist $\beta$ and $\gamma$, such that  $\beta >\gamma>0$,
and,
\begin{equation}\label{Eq_Cone-cond}
\forall \; x \in F, \; \gamma R(\|x\|) \leq \Phi(x) \leq \beta
R(\|x\|).
\end{equation}
For obvious reasons, $R$ will be called sometimes a growth function.
\end{itemize}
\end{definition}

\begin{remark}\label{Rem_mesurability}
\begin{itemize}
\item[(a)] It is interesting to note that the closure
$K:=\overline{\mathcal{E}_F^R}$, is a closed cone with non-empty
interior in the Banach space $X=C^0(F,\mathbb{R})$ of continuous
functions $\Psi:F \rightarrow \mathbb{R}$, endowed with the
compact-open topology \cite{TOPOhirch}, {\it i.e.} $K+K\subset K$,
$tK\subset K$ for every $t \geq 0$, $K\cap (-K)=\{0_X\}$ and Int
$K\neq \emptyset$.

\item[(b)] Note that the results obtained in this article could be derived
with weaker assumptions than in $(G_3)$, such as relaxing
(\ref{Eq_Cone-cond}) for $\|x\| \geq \nu $ for some $\nu>0$, and
assuming measurability on $R$ and $\Phi$ (with respect to the Borel
$\sigma$-algebras of $\mathbb{R}^+$ and $F$ respectively) instead of
continuity. However, further properties have to be derived in order
to extend appropriately the approach developed in this paper. For
instance assuming only measurability of $R$, it can be proved, since
$R$ is assumed to be subadditive, that $R$ is bounded on compact
subsets of $\mathbb{R}^+$, {\it e.g.} \cite[lemma 1,
p. 167]{TOPO_gosselin}; a property that would appear to be important for extending
the results of this article in such a context. We leave for the
interested reader these possible extensions of the results presented
hereafter.

\item[(c)] Other generalization about $R$ could be also considered,
such as $R(u+v)\leq C(R(u)+R(v)), \; \forall \; u,v\in
\mathbb{R}^{+},$ for some $C>0$, allowing {\jr for} the fact that any positive
power of a subadditive function is subadditive in that sense; but
this condition would add complications in the proof of  Theorem
\ref{Thm_fixed-point} for instance. We do not enter in all these
generalities to make the expository less technical.
\end{itemize}

\end{remark}

We need also to consider a function $r:\mathbb{R}^+ \rightarrow
\mathbb{R}^+$, verifying the following assumptions.
\begin{Aonr}
We assume that $r(x)=0$ if and only if $x=0$, $r$ is continuous at
$0$, $r$ is nondecreasing, subadditive and for some statements we
will assume {\mkk furthermore} that,
\begin{itemize}
\item[(A$_r$)] $r$ is lower semi-continuous for the {\mkk pointwise} 
convergence on $F$, {\it i.e.},
$$ {\mkk \forall \; x\in F, \;  r(\liminf_{n\rightarrow +\infty} \|f_n(x)\|)\leq \liminf_{n\rightarrow +\infty} r(\|f_n(x)\|)}, $$
for any sequence $\{f_n\}_{n\in \mathbb{N}}$ of self-mappings of
$F$.
\end{itemize}
\end{Aonr}

\begin{CrR}
Finally, we will consider the following cross condition between the
growth function $R$ and the function $r$,
\begin{itemize}
\item[(C$_{r,R}$)] $\exists \; a >0, \; \exists \; b>0, \;
\forall \; u \in \mathbb{R}^+, \; R(u)\leq ar(u) +b.$
\end{itemize}
\end{CrR}

As simple example of functions $\Phi, R$, and $r$ satisfying the
above conditions (including (A$_r$)), we can cite $r(u)=u$,
$\Phi(x)=R(\|x\|)=\sqrt{\|x\|}+1$, that will be used to illustrate
the main theorem of this article later on; see subsection
\ref{SEC_example}.

Hereafter in this subsection, condition (A$_r$) is not required. We
introduce now  the following functional on Hom($F$) with possible
infinite values,
\begin{equation}\label{Eq_def-semi-norm}
  |\cdot|_{\Phi,r}: \left\{ \begin{array}{l}
  \textrm{Hom}(F) \rightarrow \overline{\mathbb{R}^+}\\
   f\mapsto |f|_{\Phi,r}:=\underset{x\in F}\sup\Big(
   \frac{r(\|f(x)-x\|)}{\Phi(x)}\Big).
\end{array} \right.
\end{equation}
Note that,
$$|f|_{\Phi,r}=0 \mbox{ if  and only if } f=\textrm{Id}_F \mbox{ (separation condition)},$$
where Id$_F$ denotes the identity map of $F$.

\begin{definition}\label{Def_Koop}
The Koopman operator {\mkk with domain} $\mathcal{E}^R_F$, associated with $f\in$
{\rm Hom}($F$), is defined as the operator $U_f$ given by:
\begin{equation}
  U_f: \left\{ \begin{array}{l} \mathcal{E}_F^R \rightarrow
  {\mkk C^0(F,\mathbb{R}^+)}\\
\Phi \mapsto U_f(\Phi), \mbox{ where } U_f(\Phi)(x)=\Phi(f(x)), \;
\forall \; x\in F.
\end{array} \right.
\end{equation}
\end{definition}

\begin{remark}
Classically, the Koopman operator is given with other domain such as
$L^{p}(F)$ \cite{TOPOding, LM94} (generally $p=2$) and arises naturally
with the Frobenius-Perron operator in the study of ergodicity and
mixing properties of measure-preserving transformations; 
{\it e.g.} \cite{TOPO_V00,TOPOSinai_book, LM94}. The Koopman operator addresses the evolution
of phase space functions (observables), such as $\Phi$ above,
described by the linear operator $U_f$ rather than addressing a
direct study of the {\mkk nonlinear} dynamics generated by $f$. This idea has been
introduced by Koopman and von Neumann in the early 30's
\cite{TOPOKoopman_Neumann}, and has paved the road of what is called
today the spectral analysis of dynamical systems
{\mkk \cite{TOPO_V00, TOPOSinai_book, TOPO_goodson, LM94}}. We propose here more specifically to link
this spectral analysis with the {\mkk topological problem of conjugacy \cite{TOPOSmale0, TOPOSmale}, from an abstract point of view. Let us mention nevertheless, that related relationships are known to exist in certain instances, but a general treatment of the question is missing to the best of the authors' knowledge.  For example,  it is known that ergodic transformations of a compact manifold are  semiconjugate    to a  rotation on the circle if and only if there exists a non-constant eigenfunction of the Koopman operator associated with the eigenvalue $e^{-2i\pi \omega}$ for some $\omega \in \mathbb{Q}$; {\it e.g.} \cite[Proposition 8]{TOPO_Mezic1}. We consider hereafter conjugacy problems between elements of a single group of transformations (Hom($F$)), rather than semiconjugacy problems between elements of different groups \cite{TOPOkat}. Our phase space will be also always assumed to be unbounded (and thus non-compact).}
\end{remark}

{\mkk \begin{remark}\label{Rem_Koopman_stable}
Note that in general $U_f$, as defined in Definition \ref{Def_Koop}, does not leave stable  $\mathcal{E}^R_F$, since $U_f(\Phi)$ is not guaranteed to satisfy (G$_3$) for $\Phi \in \mathcal{E}^R_F$. In fact a direct analysis shows that in order to have  $U_f(\Phi)$  to satisfiy (G$_3$) it requires   restrictions on $f$ and $R$ that we want to avoid\footnote{{\mkk For instance if we assume $c\|x\|\leq \|f(x)\| \leq C \|x\|$ for all $x \in F$, and $R$ to be furthermore increasing and quasi-homogeneous \cite{TOPO_Rosen}, then $U_f(\Phi)$ satisfies (G$_3$). Assuming furthermore that $\|f(x)\|\rightarrow \infty$ as $\|x\|\rightarrow \infty$, then  $U_f(\Phi)$ satisfies (G$_2$) and $U_f$ leaves thus stable $\mathcal{E}^R_F$ in such a case.}}.  However if there exists some positive constants $c(f)$  and $C(f)$ such that $c(f) \Phi\leq U_f(\Phi) \leq C(f) \Phi$ for a particular $\Phi$ in $\mathcal{E}^R_F$, then we can conclude that $U_f(\Phi)$ lives in $\mathcal{E}^R_F$. In the proposition below, we derive when  $f \in \mathbb{H}_{\Phi,r}$ such a upper bound valid for all $\Phi$ in $\mathcal{E}^R_F$; see \eqref{Eq_Upper-bound-Koopman2}. The lower bound will be naturally verified for generalized eigenfunctions of $U_f$ in the sense of  Definition \ref{Def_Koopman_gene_eigen} below; making such eigenfunctions elements of $\mathcal{E}^R_F$.
\end{remark}}

We can now state the following proposition.

\begin{proposition}\label{Thm_group}
Consider $R$ given as in Definition \ref{Def_SpaceE_F},  and
$\Phi\in \mathcal{E}_F^R$. Let $r$ satisfy the above assumptions
except (A$_r$), and such that (C$_{r,R}$) is satisfied. Introduce
the following subset of {\rm Hom}($F$),
\begin{equation}\label{Eq_Hphir}
\mathbb{H}_{\Phi,r}:=\{f\in {\rm Hom}(F)\; : \;
|f|_{\Phi,r}<\infty, \mbox{ and } |f^{-1}|_{\Phi,r}<\infty\}.
\end{equation}
Then  ($\mathbb{H}_{\Phi,r},\circ$) is a subgroup of
({\rm Hom}($F$),$\circ$) and, for any $f \in \mathbb{H}_{\Phi,r}$, the
Koopman operator{\mkk $,U_f,$} associated with $f$ is a bounded operator on $\mathcal{E}^R_F$ which satisfies,
\begin{equation}\label{Eq_Upper-bound-Koopman2}
\forall \; \Phi \in \mathcal{E}^R_F, \; U_f(\Phi)\leq \Lambda(f)
\Phi,
\end{equation}
with,
\begin{equation}\label{Eq_defLambda}
\Lambda(f):=a\beta |f|_{\Phi,r}+bm^{-1}\beta + \beta
\gamma^{-1} \mkk{<\infty},
\end{equation} 
and where the constants appearing here are as
introduced above.
\end{proposition}
\begin{proof}

We first prove the subgroup property. Let $x$ be arbitrary in $F$,
and $f$, $g$ in $\mathbb{H}_{\Phi,r}$. Then,
\begin{equation}\label{Eq_int1}
 \begin{split}
\frac{r(\|f\circ g^{-1}(x)-x\|)}{\Phi(x)}\leq \frac{r(\|f\circ
g^{-1}(x)-g^{-1}(x)\|)}{\Phi(g^{-1}(x))}& \cdot
\frac{\Phi(g^{-1}(x))}{\Phi(x) }  \\&+\frac{r(\|
g^{-1}(x)-x\|)}{\Phi(x)}.
\end{split}
\end{equation}
From (G$_3$) and the subadditivity of $R$,
$$ \Phi(g^{-1}(x))\leq \beta (R(\|g^{-1}(x) -x\|)+R(\|x\|)),$$
and since $R(\|x\|)\leq \gamma^{-1} \Phi(x),$ we get by using
$(C_{r,R})$  and (G$_1$),
\begin{equation}\label{Eq_int3}
    \begin{split}
\frac{\Phi(g^{-1}(x))}{\Phi(x)}\leq
&\beta\left(\frac{a\;r(\|g^{-1}(x)-x)\|)}{\Phi(x)}+\frac{b}{\Phi(x)}\right)+
\beta{\gamma}^{-1}\\
\leq&  C:=a\beta |g^{-1}|_{\Phi,r}+bm^{-1}\beta+\beta{\gamma}^{-1},
\end{split}
\end{equation}
with $C$ finite since $|g^{-1}|_{\Phi,r}$ exists by definition of
$\mathbb{H}_{\Phi,r}.$

Going back to (\ref{Eq_int1}) we deduce that,
\begin{equation}\label{Eq_int5}
\frac{r(\|f\circ g^{-1}(x)-x\|)}{\Phi(x)}\leq C|f|_{\Phi,r}
+|g^{-1}|_{\Phi,r} < \infty,
\end{equation}
which concludes that $f\circ g^{-1} \in \mathbb{H}_{\Phi,r},$ and
$\mathbb{H}_{\Phi,r}$ is a subgroup of Hom($F$). The proof of
(\ref{Eq_Upper-bound-Koopman2}) consists then just in a
reinterpretation of (\ref{Eq_int3}).
\end{proof}

\begin{remark}
Fairly general homeomorphisms are encompassed by the groups,
$\mathbb{H}_{\Phi,r}$, introduced above. For instance, in the
special case $\Phi(x)=R(\|x\|):=\|x\|+1, $ and $r(x)=x$, denoting by
$\mathbb{H}_0$ the group $\mathbb{H}_{\Phi,r},$ and
$|\cdot|_{\Phi,r}$ by $|\cdot|_0$ for that particular choice of
$\Phi,$ and $r$, the following two classes of homeomorphisms belong
to $\mathbb{H}_0$ and exhibit non-trivial dynamics.

\begin{itemize}
\item[(a)] Mapping $f$ of $\mathbb{R}^d$ which
are perturbation of linear mapping in the following sense:
\begin{equation}\label{general_elt}
f(x)=Tx + \varphi(x),
\end{equation}
with $T$ a {\it linear automorphism} of $\mathbb{R}^d$ and $\varphi$
a C$^1$ map which is globally Lipschitz with Lipschitz constant,
Lip($\varphi$), satisfying Lip($\varphi$)$<
||T^{-1}||_{\mathcal{L}(\mathbb{R}^d)}^{-1}$
--- that ensures $f$ to be an homeomorphism of $\mathbb{R}^d$ from
the {\it Lipschitz inverse mapping theorem} (cf. {\it e.g.}
\cite[p. 244]{TOPOIrwin}) --- and $\|\varphi(x)\| \leq C (\|x\|+1)$
(that ensures $|f|_0<\infty$) for some positive constant; and such
that the inverse of the differential of $f$ has a uniform upper
bound $M>0$ in the operator norm, {\it i.e.},
$\|[Df(u)]^{-1}\|_{\mathcal{L}(\mathbb{R}^d)}\leq M$ for every $u\in
\mathbb{R}^d$, which ensures $|f^{-1}|_0 < \infty$ by the mean value
theorem. For instance, $f(x)=x+\frac{1}{2}\log(1+x^2)$ provides such
an homeomorphism of $\mathbb{R}$. Note that every $\varphi$ that is
a C$^1$ map of $\mathbb{R}^d$ with compact support and with
appropriate control on its differential leads to an homeomorphism
 of type (\ref{general_elt}) that belongs to $\mathbb{H}_0$.

\item[(b)] Extensions of the preceding examples to the class of non-smooth (non C$^1$) perturbations of linear
automorphisms can be considered; exhibiting non-trivial
homeomorphisms of $\mathbb{H}_0$. For instance, the two-parameter
family of homeomorphisms, $\{\mathbb{L}_{a,b} \;, \;
b\in\mathbb{R}\backslash\{0\}, \; a \in \mathbb{R}\}$, known as the
Lozi maps family \cite{TOPOLozi}:
\begin{equation}\label{Lozi}
\left.
\begin{array}{l}
\mathbb{L}_{a,b}: \mathbb{R}^2 \longrightarrow \mathbb{R}^2\\
 \; \; \; \quad \; (x,y) \mapsto (1-a|x|+y, bx)
\end{array} \right.,
\end{equation}
where $|\cdot |$ denotes the absolute value here; constitutes a
family of elements of $\mathbb{H}_{0}$. Indeed, it is not difficult
to show that $\mathbb{L}_{a,b}$ and its inverse, for $b\neq 0$,
$\mathbb{M}_{a,b}:(u,v)\mapsto (\frac{1}{b}v,
-1+u+\frac{a}{|b|}|v|)$ have finite $|\cdot|_0$-values. This family
shares similar properties with the H\'enon maps family. For instance
there exists an open set in the parameter space for which
generalized hyperbolic attractors {\jr exist}  \cite{TOPOMisiu}.
\end{itemize}
\end{remark}

We introduce now the following functional on
$\mathbb{H}_{\Phi,r}\times\mathbb{H}_{\Phi,r},$

\begin{equation}\label{Eq_def_rho}
\rho_{\Phi,r}(f,g):=\max(|f\circ g^{-1}|_{\Phi,r}, |f^{-1}\circ
g|_{\Phi,r}),
\end{equation}
which is well-defined by Proposition \ref{Thm_group} and
non-symmetric. Since obviously $\rho_{\Phi,r}(f,g)\\ \geq 0$ whatever
$f,$ and $g$, and $\rho_{\Phi,r}(f,g)=0$ if and only if $f=g,$ then
$\rho_{\Phi,r}$ is in fact a premetric on $\mathbb{H}_{\Phi,r}$.
Note that hereafter, we will simply denotes $f\circ g$ by $fg$. Due
to the non-symmetric property, two natural {\jr types} of ``balls'' can be
defined with respect to the premetric $\rho_{\Phi,r}$. More
precisely:

\begin{definition}
An open $\rho_{\Phi,r}$-ball of center $f$ to the right (resp. left)
and radius $\alpha >0$ is the subset of $\mathbb{H}_{\Phi,r}$
defined by $B^{+}_{\rho_{\Phi,r}}(f,\alpha):=\{ g\in
\mathbb{H}_{\Phi,r} \; : \; \rho_{\Phi,r}(g,f) < \alpha \}$ (resp.
$B^{-}_{\rho_{\Phi,r}}(f,\alpha):=\{ g\in \mathbb{H}_{\Phi,r} \; :
\; \rho_{\Phi,r}(f,g) < \alpha \}$).
\end{definition}

\begin{proposition}\label{Thm_Main-propertiesI}
Consider $R$ given as in Definition \ref{Def_SpaceE_F},  and
$\Phi\in \mathcal{E}_F^R$. Let $r$ satisfy the above assumptions
except (A$_r$), and such that (C$_{r,R}$) is satisfied. Then, the
premetric as defined in \eqref{Eq_def_rho} satisfies the following
properties.

\begin{itemize}
\item[(i)] For every $f,g, h,$ in $\mathbb{H}_{\Phi,r}$,  the
following relaxed triangle inequality holds,
\begin{equation}\label{Eq_relaxed_ineq2}
\rho_{\Phi,r}(f,g)\leq a \beta \rho_{\Phi,r}(f,h) \rho_{\Phi,r}(h,g)
+(bm^{-1}\beta +\beta \gamma^{-1}) \rho_{\Phi,r}(f,h)+
\rho_{\Phi,r}(h,g).
\end{equation}

\item[(ii)] The following families of subsets of
$\mathbb{H}_{\Phi,r}$,
$$\mathfrak{T}^{+}(\rho_{\Phi,r}):=\{\mathcal{H} \subset \mathbb{H}_{\Phi,r} \;
 : \; \forall \; f\in \mathcal{H}, \exists \; \alpha >0, \; B^{+}_{\rho_{\Phi,r}}(f,\alpha) \subset \mathcal{H} \} $$
and,
$$ \mathfrak{T}^{-}(\rho_{\Phi,r}):=\{\mathcal{H} \subset \mathbb{H}_{\Phi,r} \;
 : \; \forall \; f\in \mathcal{H}, \exists \; \alpha >0, \; B^{-}_{\rho_{\Phi,r}}(f,\alpha) \subset \mathcal{H} \} $$
are two topologies on $\mathbb{H}_{\Phi,r}$.

\item[(iii)] For all $f\in \mathbb{H}_{\Phi,r}$, for all $\alpha^*
>0$, and for all $g\in B^{-}_{\rho_{\Phi,r}}(f,\alpha^*)$, the
following property holds:
$$\Big(\rho_{\Phi,r}(f,g)< \frac{\alpha^{*}}{bm^{-1} \beta + \beta \gamma^{-1}})\Big)\Rightarrow (\exists \; \alpha>0, \;
B^{-}_{\rho_{\Phi,r}}(g,\alpha) \subset
B^{-}_{\rho_{\Phi,r}}(f,\alpha^*)), $$ and thus for all $f\in
\mathbb{H}_{\Phi,r}$, $\underset{\alpha>0}\bigcup
\;B^{-}_{\rho_{\Phi,r}}(f,\alpha)$ is a fundamental system of
neighborhoods of $f$, which renders
$\mathfrak{T}^{-}(\rho_{\Phi,r})$ first-countable. {\jr An analogous}
statement holds with ``+'' instead of ``-''.

\item[(iv)] Let {\jr $\overline{\mathbb{H}}_{\Phi,r}^{(-)}$} denote the closure of $\mathbb{H}_{\Phi,r}$ for the topology
$\mathfrak{T}^{-}(\rho_{\Phi,r})$, then
$${\jr \overline{\mathbb{H}}_{\Phi,r}^{(-)}}\bigcap {\rm Hom}(F)\subset
\mathbb{H}_{\Phi,r}.$$

\end{itemize}
\end{proposition}

\begin{remark}
Proof of (iii) below shows that an arbitrary open
$\rho_{\Phi,r}-$ball (centered to the right or left) is not
necessarily open in the sense of not being an element $\mathcal{H}$
of $\mathfrak{T}^{+/-}(\rho_{\Phi,r})$, since $b>0$ and $\gamma <
\beta$.
\end{remark}

\begin{proof}
We first prove (i). Using the triangle inequality for $\|\cdot \|$
and subadditivity of $r$, it is easy to note that for all $x\in F$, {\mkk and  all $f,g,h\in \mathbb{H}_{\Phi,r}$,}

\begin{equation}\label{Eq_int6}
\frac{r(\|fg^{-1}(x)-x\|)}{\Phi(x)} \leq
\frac{r(\|fh^{-1}hg^{-1}(x)-hg^{-1}(x)\|)}{\Phi(x)} +
|hg^{-1}|_{\Phi,r}.
\end{equation}

From the following trivial equality,
$$r(\|fh^{-1}hg^{-1}(x)-hg^{-1}(x)\|)=r(\|fh^{-1}hg^{-1}(x)-hg^{-1}(x)\|)\frac{\Phi(hg^{-1}(x))}{\Phi(hg^{-1}(x))},$$
we deduce from Proposition \ref{Thm_group}, that,

$$\frac{r(\|fg^{-1}(x)-x\|)}{\Phi(x)} \leq \Lambda(hg^{-1}) |fh^{-1}|_{\Phi,r}, $$
{\mkk where $\Lambda(hg^{-1})$ is well-defined since $ \mathbb{H}_{\Phi,r}$ is a subgroup.} {\mkk This last inequality} reported in (\ref{Eq_int6}) gives {\mkk then},
\begin{equation}\label{Eq_int7}
\underset{x\in F}\sup \Big(\frac{r(\|fg^{-1}(x)-x\|)}{\Phi(x)}\Big)
\leq \Lambda(hg^{-1}) |fh^{-1}|_{\Phi,r} + |hg^{-1}|_{\Phi,r},
\end{equation}
leading to (\ref{Eq_relaxed_ineq2}), by re-writing appropriately
(\ref{Eq_int7}) and repeating the computations with the
substitutions $f\leftarrow f^{-1}$, $g^{-1} \leftarrow g$ and
$h^{-1} \leftarrow h$ for the estimation of $|f^{-1}g|_{\Phi,r}$.

The proof of (ii) is just a classical ``game" with the axioms of a
topology and is left to the reader.

We prove now (iii), only for $\mathfrak{T}^{-}(\rho_{\Phi,r})$; the
proof for $\mathfrak{T}^{+}(\rho_{\Phi,r})$ being  a repetition. Let
$f \in \mathbb{H}_{\Phi,r}$ and $\alpha^{*}>0$. Let $g\in
B_{\rho_{\Phi,r}}^{-}(f,\alpha^*),$ then from
\eqref{Eq_relaxed_ineq2} we get for all $h\in \mathbb{H}_{\Phi,r}$,
\begin{equation}\label{Eq_int8}
\rho_{\Phi,r}(f,h) \leq a\beta \rho_{\Phi,r}(f,g)\rho_{\Phi,r}(g,h)
+ (bm^{-1}\beta +\beta \gamma^{-1})\rho_{\Phi,r}(f,g)+
\rho_{\Phi,r}(g,h).
\end{equation}

We seek now the existence of $\alpha>0$ such that
$B_{\rho_{\Phi,r}}^{-}(g,\alpha)\subset
B_{\rho_{\Phi,r}}^{-}(f,\alpha^*)$. Denoting $\rho_{\Phi,r}(f,g)$ by
$\alpha'$, such a problem of existence is then reduced from Eq.
\eqref{Eq_int8} to the existence of a solution $\alpha>0$ of,
\begin{equation}\label{Eq_system1}
a\beta  \alpha \alpha' + (bm^{-1}\beta +\beta \gamma^{-1})\alpha' + \alpha < \alpha^{*}.\\
\end{equation}
A necessary condition of existence is,
\begin{equation}
\alpha' < \alpha^{**}:=\frac{\alpha^*}{bm^{-1}\beta +\beta
\gamma^{-1}}
\end{equation}
that turns out to be sufficient since any $\alpha >0$ satisfying,
\begin{equation}
\alpha < \frac{\alpha^*-(bm^{-1}\beta +\beta \gamma^{-1})
\alpha'}{1+a\beta \alpha'},
\end{equation}
is a solution because the RHS is positive.

The second part of (iii) is a reinterpretation of the result just
obtained. Indeed, we have  proved that for all $f \in
\mathbb{H}_{\Phi,r},$ and for all $\alpha^*>0$ there exists $ 0 <
\alpha^{**}< \alpha^*$ (since $\gamma < \beta$ and $b > 0$ by
definition), such that $B_{\rho_{\Phi,r}}^{-}(f,\alpha^{**}) \in
\mathfrak{T}^{-}(\rho_{\Phi,r})$.

Now if we introduce $\mathbb{B}(f):=\underset{\alpha^*>0}\bigcup
B_{\rho_{\Phi,r}}^{-}(f,\alpha^{*}),$ then the family
$\mathcal{F}(f)$ of subsets of $\mathbb{H}_{\Phi,r}$ defined by,
$$\mathcal{F}(f):=\{V \in 2^{\mathbb{H}_{\Phi,r}} \; : \; \exists \; B \in \mathbb{B}(f), \mbox{ s. t. } V\supset B \},$$
is a family of neighborhoods of $f$, since every $V \in
\mathcal{F}(f)$ contains by definition a subset
$B_{\rho_{\Phi,r}}^{-}(f,\alpha^{*})$ (that is not necessarily
open), and therefore a subset of type
$B_{\rho_{\Phi,r}}^{-}(f,\alpha^{**})$ which is open from what
precedes. {\jr Thus first countability  naturally holds}  and the proof
of (iii) is complete.

We prove now (iv). Let $f \in
{\jr \overline{\mathbb{H}}_{\Phi,r}^{(-)}}\bigcap \textrm{Hom}(F),$ then
by the property (iii), there exists a sequence $\{ f_n\}_{n\in
\mathbb{N}} \in \mathbb{H}_{\Phi,r}^{\mathbb{N}}$, such that
$\rho_{\Phi,r}(f,f_n)\underset{n\rightarrow \infty}\rightarrow 0.$
Then by definition of $\rho_{\Phi,r}$, we get in particular that
$|ff_n^{-1}|_{\Phi,r}$ and $|f^{-1}f_n|_{\Phi,r}$ exist from which
we deduce  that $|f|_{\Phi,r}$ and $|f^{-1}|_{\Phi,r}$ exist, since
$\mathbb{H}_{\Phi,r}$ is a group.

\end{proof}

\subsection{Closure properties for extended premetric on Hom($F$)}\label{Sec_prelim2}
In this section, we extend the closure property (iv) of Proposition
\ref{Thm_Main-propertiesI} to Hom($F$) itself (cf. Proposition
\ref{Thm_Closure1}) and prove a cornerstone proposition (Proposition
\ref{Thm_Closure2}) concerning the convergence in
$\mathfrak{T}^-(\rho'_{\Phi,r})$ of sequences taking values in
Hom($F$), where $\rho'_{\Phi,r}$ denotes the extension of the
premetric $\rho_{\Phi,r}$ to Hom($F$)$\times$Hom($F$).

The result described in Proposition \ref{Thm_Closure1} will allow us
to make precise conditions for which the solution of the fixed point
theorem proved in the next section, lives in $\mathbb{H}_{\Phi,r}$.
This specific result is not fundamental for the proof of the main
theorem of conjugacy of this article, Theorem \ref{THM_conjug};
whereas Proposition \ref{Thm_Closure2} will play an essential role
in the proof of the fixed point theorem, Theorem
\ref{Thm_fixed-point}, and by the  way in the proof of Theorem
\ref{THM_conjug}. Important related concepts such as the one of
incrementally bounded sequence are also introduced in this
subsection.

We define $\rho'_{\Phi,r}$  as the extension of the premetric
$\rho_{\Phi,r}$ to Hom($F$)$\times$Hom($F$), by classically allowing
$\rho'_{\Phi,r}$ to take values in the extended real line
$\overline{\mathbb{R}}:=\mathbb{R}\cup \{\infty\}$ instead of
$\mathbb{R}$; with the  usual extensions of the arithmetic
operations\footnote{Note that similarly, $|\cdot|_{\Phi,r}$ may be
extended in such a way, but it is important to have in mind that
$|f|'_{\Phi,r}=\infty$ and $|g|'_{\Phi,r}=\infty$ do not necessarily
imply that $\rho'_{\Phi,r}(f,g)=\infty$. {\jr This is, for instance,} the case
for $f(x)=g(x)=Ax$, with $A\in Gl_{d}(\mathbb{R}), A\neq I_{d}$,
$r(x)=x$, and $\Phi(x)=\sqrt{\|x\|}+1$.}.

According to this basic extension procedure, it can be shown that
$\mathfrak{T}^{+/-}(\rho'_{\Phi,r})$ is a topology on Hom($F$) and
that Proposition \ref{Thm_Main-propertiesI} can be reformulated for
$\rho'_{\Phi,r}$ with the appropriate modifications. Note that{\jr ,} in
particular{\jr ,} the relaxed triangle inequality \eqref{Eq_relaxed_ineq2}
holds for $\rho'_{\Phi,r}$ and any $f,g$ and $h$ in Hom($F$).
{\mkk Indeed, either $fg^{-1}$ and $f^{-1}g$ {\jr both} belong to
$\mathbb{H}_{\Phi,r}$ {\jr in which case} \eqref{Eq_relaxed_ineq2} obviously holds;
{\jr or}  $fg^{-1}$ {\jr and}  $ f^{-1}g$ {\jr do not both} belong to
$\mathbb{H}_{\Phi,r}$.  {\jr In the latter case,} because of the group structure of $\mathbb{H}_{\Phi,r}$,  at least one element {\jr in} $\{fh^{-1}, f^{-1}h, hg^{-1}, h^{-1}g\}$ does
not belong to $\mathbb{H}_{\Phi,r}$. {\jr This leads to the conclusion that the inequalities}  \eqref{Eq_relaxed_ineq2} {\jr still hold} for
$\rho'_{\Phi,r}$.}

We are now in {\jr a} position to introduce contingent conditions to our
framework that are required to obtain closure type results in
Hom($F$). These conditions {\mk possess} the particularity to hold in
$\mathfrak{T}^{+}(\rho'_{\Phi,r})$ for any sequences involved in the
closure problem related to the topology
$\mathfrak{T}^{-}(\rho'_{\Phi,r})$. More precisely we introduce the
following concepts.

\begin{definition}\label{Def_closure1}
We say that a sequence $\{f_n\}$ of elements of {\rm Hom}($F$) is
incrementally bounded in $\mathfrak{T}^{+}(\rho'_{\Phi,r})$ with
respect to $q$, for some $q\in \mathbb{N}$, if and only if:
$$\exists \; C^{+}_q, \; \forall \; p\in \mathbb{N}, \; (p\geq q) \Rightarrow (\rho'_{\Phi,r}(f_p,f_q)\leq C^{+}_q).$$

When  {\jr no further conditions on $q$ are assumed, we say}  that the sequence is
incrementally bounded.

Furthermore, the sequence is said to be uniformly bounded in
$\mathfrak{T}^{+}(\rho'_{\Phi,r})$ if and only if:
$$\exists \; C^{+}, \; \forall \; p,q\in \mathbb{N}, \; (p\geq q) \Rightarrow (\rho'_{\Phi,r}(f_p,f_q)\leq C^{+}_q).$$
The ``(-)-statements'' consist of changing the role of $p$ and $q$ in
the above statements. We denote by $\mathcal{IB}^{+}$ the set of
incrementally bounded sequences in
$\mathfrak{T}^{+}(\rho'_{\Phi,r})$ and by $\mathcal{IB}^{+}_u$ its
subset constituted only by uniformly incrementally bounded
sequences.

\end{definition}

{\mkk \begin{definition}\label{Def_closure2}
We define $\overline{{\rm Hom}(F)}^{(-),b^{+}}$ to be the set {\jr consisting of all the limit points  in $\mathfrak{T}^{-}(\rho'_{\Phi,r})$  of
sequences $\{f_n\}\in$ {\rm Hom}($F$),} 
such that,
$$\exists \; n_0 \in \mathbb{N}\; :\; f_{n_0} \in \mathbb{H}_{\Phi,r}\;, $$
for which $\{f_n\}$ is incrementally bounded in
$\mathfrak{T}^{+}(\rho'_{\Phi,r})$ with respect to $n_0$.
\end{definition}}

We have then the {\mkk following} important proposition that completely characterizes
the limits in $\mathfrak{T}^{-}(\rho'_{\Phi,r})$ (that leave
Hom($F$) stable) of sequences of {\mkk elements of} Hom($F$) which are incrementally
bounded in $\mathfrak{T}^{+}(\rho'_{\Phi,r})$ with respect to some
$n_0$ for which $f_{n_0}$ belongs to $\mathbb{H}_{\Phi,r}$.
\begin{proposition}\label{Thm_Closure1}
Let $\overline{{\rm Hom}(F)}^{(-),b^{+}}$ be as introduced in Definition
\ref{Def_closure2}, then,

\begin{equation}
\overline{{\rm Hom}(F)}^{(-),b^{+}} \bigcap {\rm Hom}(F) \subset
\mathbb{H}_{\Phi,r}.
\end{equation}
\end{proposition}

\begin{proof}
Let $f\in \overline{{\rm Hom}(F)}^{(-),b^{+}} \bigcap {\rm Hom}(F)$, we want to
show that $|f|_{\Phi,r}$
 and $|f^{-1}|_{\Phi,r}$ exist. By assumptions, there exists
 $\{f_n\} \in \mathcal{IB}^{+}$, such that
 $\rho'_{\Phi,r}(f,f_n)\underset{n\rightarrow \infty}\longrightarrow
 0$. Consider $n_0$ such that $f_{n_0} \in \mathbb{H}_{\Phi,r}$
 resulting from the definition of $\overline{{\rm Hom}(F)}^{(-),b^{+}}$.
 From \eqref{Eq_relaxed_ineq2},
\begin{equation}\label{Eq_int9}
\begin{split}
 \rho'_{\Phi,r}(f_p, Id_F) \leq a\beta& \rho'_{\Phi,r}(f_p, f_{n_0})\rho'_{\Phi,r}(f_{n_0},
 Id_F)\\
 &+(bm^{-1}\beta +\beta \gamma^{-1})\rho'_{\Phi,r}(f_p, f_{n_0}) + \rho'_{\Phi,r}(f_{n_0},
 Id_F).
 \end{split}
 \end{equation}

Since $\{f_n\}$ is incrementally bounded in
$\mathfrak{T}^{+}(\rho'_{\Phi,r})$ with respect to $n_0$, then from
\eqref{Eq_int9},
 we deduce that the real-valued sequence $\{|f_{n}|_{\Phi,r}
\}_{n\geq n_0}$ is bounded.

Besides, for any $x\in F$, and any $n\geq n_0$,
$$r(\|f(x)-x\|) \leq |ff^{-1}_n|_{\Phi,r} \Phi(f_n(x)) +  r(\|f_n(x)-x\|),$$
and therefore by using the estimate \eqref{Eq_Upper-bound-Koopman2}
in Proposition \ref{Thm_group},

\begin{equation}\label{Eq_int11}
\frac{r(\|f(x)-x\|)}{\Phi(x)} \leq |ff^{-1}_n|_{\Phi,r}\cdot
\Lambda(|f_n|_{\Phi,r}) + |f_n|_{\Phi,r},
 \end{equation}
since $\Lambda(|f_n|_{\Phi,r})$ is well defined for $n\geq n_0$
because $|f_n|_{\Phi,r}$ exist for such $n$. We get then trivially
that $\{\Lambda(|f_n|_{\Phi,r})\}_{n\geq n_0}$ is bounded because
$\{|f_n|_{\Phi,r}\}_{n\geq n_0}$ is. Now since $|ff_n^{-1}|_{\Phi,r}
\underset{n\rightarrow \infty}\longrightarrow 0$ by assumption, we
then deduce that $|f|_{\Phi,r}$ {\jr exists}  by taking $n$ sufficiently
large in \eqref{Eq_int11}. To conclude,  it suffices to note that
\eqref{Eq_int9}
shows thanks to the incrementally bounded property, that
$\{|f_n^{-1}|_{\Phi,r}\}_{n\geq n_0}$ is bounded as well, leading to
the boundedness of $\{\Lambda(|f_n^{-1}|_{\Phi,r})\}_{n\geq n_0}$
which by repeating similar estimates leads to the existence of
$|f^{-1}|_{\Phi,r}$ by using the assumption: $|f^{-1}f_n|_{\Phi,r}
\underset{n\rightarrow \infty}\longrightarrow 0$.
\end{proof}

In the sequel we will need sometimes to use the following property
verified by the function $r$. Let $r$ satisfy  the  condition{\mkk s} of the
preceding subsection and let $G$ denote a continuous function
$G:F\rightarrow E$, then for any $K$ compact {\mk subset} of $F$, there exists
$x_K\in K, $ such that $r\Big( \underset{x\in K}\sup
\|G(x)\|\Big)=r(\|G(x_K)\|)=\underset{x\in K}\sup \; r(\|G(x)\|),$
since $r$ is increasing. Since we will need this property of $r$
later, we make it precise as the condition,
\begin{equation*}
(S)\;\;:\;\mbox{for all compact set } K\subset F, \;
r\Big(\underset{x\in K}\sup \|G(x)\|\Big) = \underset{x\in K}\sup \;
r(\|G(x)\|),
\end{equation*}
for every continuous function $G:F\rightarrow E$;  {\mkk condition} which holds
therefore for $r$ as defined above. In what follows, $\rho_{\Phi,r}$ will refer for both $\rho_{\Phi,r}$
when applied to elements of $\mathbb{H}_{\Phi,r}$, and to
$\rho'_{\Phi,r}$ when applied to elements of Hom($F$), without any
sort of confusion.

 Let us now introduce the following concept of Cauchy sequence in
$\mathbb{H}_{\Phi,r}$ or more generally in Hom($F$), which is
adapted to our framework.

\begin{definition}\label{Def_Cauchy}
A sequence $\{f_n\}$ in $\mathbb{H}_{\Phi,r}$ or in ${\rm Hom}(F)$, is
called $\rho_{\Phi,r}^{+}$-Cauchy (resp.
$\rho_{\Phi,r}^{-}$-Cauchy), if the following condition holds,
$$\forall \; \epsilon >0, \; \exists \; N \in \mathbb{N}, \; (p\geq q \geq N) \Rightarrow (\rho_{\Phi,r}(f_p,f_q)\leq \epsilon). $$
$$(\mbox{resp.} \; \forall \; \epsilon >0, \; \exists \; N \in \mathbb{N}, \; (q\geq p \geq N) \Rightarrow (\rho_{\Phi,r}(f_p,f_q)\leq \epsilon). $$
\end{definition}

\begin{remark}
Note that, since $\rho_{\Phi,r}$ is not symmetric, the role of $p$
and $q$ are not symmetric as well, to the contrary of the classical
definition of a Cauchy sequence in a metric space.

\end{remark}

\begin{remark}\label{Rem_Cauchy_is_IBu}
By definition, every sequence $\{f_n\}$ in $\mathbb{H}_{\Phi,r}$
which is $\rho_{\Phi,r}^{+}-$Cauchy (resp.
$\rho_{\Phi,r}^{-}-$Cauchy) belongs to $\mathcal{IB}_u^{+}$ (resp.
$\mathcal{IB}_u^{-}$). However, a sequence $\{f_n\}$ in Hom($F$)
which is $\rho_{\Phi,r}^{+}-$Cauchy is not an element of
$\mathcal{IB}_u^{+}$ in general but is an element of
$\mathcal{IB}^{+}.$

\end{remark}

We are now in  position to prove the following cornerstone
proposition concerning the convergence in
$\mathfrak{T}^{-}(\rho_{\Phi,r})$ of $\rho_{\Phi,r}^{+}$-Cauchy
sequences in Hom($F$). {\mk We recall that a topological space is defined to be $\sigma$-compact if and only if it is the union of a 
countable family of compact subsets \cite{Kelley}.}

\begin{proposition}\label{Thm_Closure2}
Assume that $F$ is an {\mk unbounded} subset of $E$ which is locally
connected, $\sigma$-compact and locally compact\footnote{Note that,
since $E$ is a finite dimensional normed vector space, which is
locally compact Haussdorf space, if $F$ is an open or closed subset
of $E$, it is locally compact; cf. {\it e.g} 
\cite[\S 3.18.4, p. 66]{TOPODieud}.}. Consider $R$ given as in Definition
\ref{Def_SpaceE_F},  and $\Phi\in \mathcal{E}_F^R$. Let $r$ satisfy
the above assumptions including (A$_r$), and such that (C$_{r,R}$)
is satisfied. Let $\{f_n\}$ be a sequence in ${\rm Hom}(F)$. If the
following conditions hold:
\begin{itemize}
\item[(C$_1$)] For every compact $K\subset F$, the sequence of the restriction of $f_n$ to $K$, $\{f_n|_K\},$ is
 bounded. {\jr The same holds for the sequence of restrictions of the inverses,} 
 $\{f_n^{-1}|_K\}$,
\item[(C$_2$)] $\{f_n\}$ is $\rho_{\Phi,r}^+$-Cauchy,
\end{itemize}
then $\{f_n\}$ converges in $\mathfrak{T}^{-}(\rho_{\Phi,r})$
towards an element of ${\rm Hom}(F)$.

If furthermore, either $\{f_n\}$ is a sequence of homeomorphisms
living in $\mathbb{H}_{\Phi,r}$, or more generally $\{f_n\}$ is such
that $f_{n_0}\in \mathbb{H}_{\Phi,r}$ for some $n_0$, then $\{f_n\}$
converges in $\mathfrak{T}^{-}(\rho_{\Phi,r})$ towards an element of
$\mathbb{H}_{\Phi,r}$.
\end{proposition}

\begin{proof}
The proof is divided in several steps.\\
{\bf Step 1.} Let $\{f_n\}$ be a sequence of homeomorphisms of $F$
fulfilling the conditions of the proposition. Let $\epsilon >0$ be
fixed. Then from (C$_2$), there exists an integer $N$ such that
$p\geq q \geq N$ implies,
\begin{equation}\label{Eq_int25}
\begin{split}
\forall \; x\in F, \;
r(\|f_p(x)-f_q(x)\|)=&r(\|f_pf_q^{-1}f_q(x)-f_q(x)\|)\\
\leq& \rho_{\Phi,r}(f_p,f_q) \Phi(f_q(x)) \leq \epsilon
\Phi(f_q(x)),
\end{split}
\end{equation}
which for every compact $K\subset F$, leads to,
\begin{equation}\label{Eq_int12}
\exists \; M_K >0, \; \forall \; x\in K, \;r(\|f_p(x)-f_q(x)\|)\leq
\epsilon M_K,
\end{equation}
by assumption (C$_1$) and the continuity of $\Phi$. Similarly,
$p\geq q \geq N$, implies,
\begin{equation}
\begin{split}
\forall \; x\in F, \;
r(\|f_p^{-1}(x)-f_q^{-1}(x)\|)=&r(\|f_p^{-1}f_qf_q^{-1}(x)-f_q^{-1}(x)\|)\\
\leq&\rho_{\Phi,r}(f_p,f_q) \Phi(f_q^{-1}(x))
 \leq \epsilon \Phi(f_q^{-1}(x)),
\end{split}
\end{equation} which, for every compact $K\subset F$, can be {\jr summarized} with
\eqref{Eq_int12}, as,
\begin{equation}\label{Eq_int13}
\begin{split}
& \exists\; M_K >0, \; r\Big(\underset{x\in
K}\sup\|f_p(x)-f_q(x)\|\Big) \leq \epsilon M_K,  \mbox{ and, } \\
& r\Big(\underset{x\in K}\sup\|f_p^{-1}(x)-f_q^{-1}(x)\|\Big)\leq
\epsilon M_K,
\end{split}
\end{equation}
by using (S) and labeling still by $M_K$ the greater constant.

 Recall  from assumptions made on the function $r$ in subsection \ref{Sec_prellim1}, that the continuity of $r$ holds at $0$ and that $r(x)=0$ if and only if
$x=0$. From that, we deduce that the sequences $\{f_n\}$ and $\{f_n^{-1}\}$
converge uniformly on each compact {\mk subset of F} towards respectively a map
$f:F\rightarrow F$ and a map $g:F\rightarrow F$. Note that by
choosing appropriately a family of compact subsets of $F$, covering
$F$, we get furthermore that $f$ and $g$ can be chosen
continuous on $F$.\\

{\bf Step 2.} Our main objective here is to show that $f=g^{-1}, $
{\it i.e.} $f\in$ Hom($F$). Since $F$ is assumed to be
$\sigma$-compact and locally compact, there exists an exhaustive
sequence of compacts sets $\{K_k\}_{k\in \mathbb{N}}$ of $F$; 
 {\it e.g.} \cite[Corollary 2.77]{TOPOAli}. From step 1,
\eqref{Eq_int13} is valid for any $p\geq q\geq N, $ with $N$ which
does not depend on the compact $K,$ and therefore we get that
$\{f_n\}$ is Cauchy for the metric $\Delta$ on Hom($F$), given by,
$$\Delta(\phi,\psi)=\delta(\phi,\psi)+\delta(\phi^{-1},\psi^{-1}),$$ where
$\delta(\phi,\psi):=\sum_{k=0}^{\infty} 2^{-k} \| \phi-\psi\|_k
(1+\| \phi-\psi\|_k)^{-1}$, and $\| \phi-\psi\|_k:=\underset{x\in
K_k}\max\|\phi(x)-\psi(x)\|$ for any compact $K_k$, and any
homeomorphisms of $F$, $\phi$ and $\psi$.

Indeed,  it suffices to note that for a given $\epsilon'>0$, from
\eqref{Eq_int13} there exists $l$, $\epsilon$ and $N_{\epsilon}$
such that,
$$\sum_{k=l+1}^{\infty}
2^{-k}  \leq \frac{\epsilon'}{4}, \mbox{ and }  \sum_{k=0}^{l}
2^{-k} \| f_p-f_q\|_k \leq \frac{\epsilon'}{4},$$ which leads to
$\delta(f_p,f_q)\leq \epsilon'/2,$ and similarly to
$\delta(f_p^{-1},f_q^{-1})\leq \epsilon'/2,$ for any $p \geq q\geq
N_{\epsilon}$.

Now, since $F$ is locally compact and locally connected from a
famous result of Arens   \cite[Theorem 4]{TOPOArens},
Hom($F$) is complete for $\Delta$ which is a metric compatible with
the compact-open topology \cite{TOPOhirch}, making Hom($F$) a Polish
group\footnote{Note that this reasoning is independent from the
choice of the metric rendering Hom($F$) complete, since it is known
that Hom($F$) has a unique Polish group structure (up to
isomorphism); see \cite{TOPOKallman}.}. Therefore $\{f_n\}$ converges
in the compact-open topology towards an element $\mathfrak{h}\in$
Hom($F$). By recalling that the compact-open topology is here
equivalent to the topology of compact convergence \cite{TOPORen}, we
obtain by uniqueness of the limit that $f=\mathfrak{h}\in$ Hom($F$);
$f$ being the limit of $\{f_n\}$ in the topology of compact
convergence from step 1.

{\bf Step 3.} Let us summarize what has been proved. We have shown
under the assumptions $(C_1)$ and $(C_2)$, that we can {\jr produce}
from any sequence $\{f_n\}$ which is $\rho_{\Phi,r}^{+}$-Cauchy, an
element $f\in {\rm Hom}(F)$, such that $\{f_n\}$ converges uniformly
to it on each compact subset of $F$, and $\{f_n^{-1}\}$ does the
same towards $f^{-1}$. In fact we can say more with respect to the
topology of convergence.

Indeed, going back to (\ref{Eq_int25}), we have that,

$$\forall \; \epsilon>0, \; \exists N \in \mathbb{N}, \; (p\geq q \geq N)\Rightarrow
\Big(\forall x\in F, \; r(\|f_p(x)-f_q(x)\|)\leq  \epsilon \cdot
\Phi(f_q(x))\Big)$$ {\it i.e.} by  making $p\rightarrow +\infty$ and
using $(A_r)$,
$$\forall \; \epsilon>0, \; \exists N \in \mathbb{N}, \; (q \geq N) \Rightarrow
\Big(\forall x\in F, \; r(\|f (x)-f_q(x)\|)\leq  \epsilon \cdot
\Phi(f_q(x))\Big),$$ which leads to,
$$ \forall \; \epsilon>0, \; \exists N \in \mathbb{N}, \; (q \geq N) \Rightarrow  (|ff_q^{-1}|_{\Phi,r} \leq \epsilon).$$
From similar estimates, we get,
$$\forall \; \epsilon>0, \; \exists N \in \mathbb{N}, \; (q \geq N) \Rightarrow  (|f^{-1}f_q|_{\Phi,r} \leq \epsilon).$$
We have thus shown the convergence of $\{f_n\}$ in
$\mathfrak{T}^{-}(\rho_{\Phi,r}),$ that is,

$$ \rho_{\Phi,r}(f,f_q)\underset{q\rightarrow +\infty }\longrightarrow 0.$$
At this stage, we have proved that $\{f_n\}$ converges in
$\mathfrak{T}^{-}(\rho_{\Phi,r})$ towards an element of Hom($F$).

{\bf Step 4.} This last step is devoted to the proof of the last
statement of the theorem concerning the membership of the limit of
$\{f_n\}$ to $\mathbb{H}_{\Phi,r}$. This fact is simply a
consequence of Proposition \ref{Thm_Main-propertiesI}-(iv) in the
case $\{f_n\}\in {\mkk (\mathbb{H}_{\Phi,r})^{\mathbb{N}}}$, and a
consequence of Remark \ref{Rem_Cauchy_is_IBu}, and Proposition
\ref{Thm_Closure1} in the case $\{f_n\}\in
({\rm Hom}(F))^{\mathbb{N}}$,  which gives in all the cases that the
limit in $\mathfrak{T}^{-}(\rho_{\Phi,r})$ of $\{f_n\}$ lives in
$\mathbb{H}_{\Phi,r}.$ The proof is therefore complete.
\end{proof}

Lastly, it is worth to note that it is only in step 3 of the
preceding proof  that was needed assumption (A$_r$), but since
Proposition \ref{Thm_Closure2} will be used in the proof of Theorems
\ref{Thm_fixed-point} and \ref{THM_conjug}, we will make a
systematic use of this assumption in the sequel.

\section{A fixed point theorem in the homeomorphisms
group}\label{Sec_Fixed_Pt_Thm} In this section we state and prove a
new fixed point theorem valid for self-mappings acting on Hom($F$),
which holds within the functional framework developed in the
preceding section. This fixed point theorem uses a contraction
mapping argument that involves a Picard scheme that has to be
controlled appropriately due to the relaxed inequality
\eqref{Eq_relaxed_ineq2}. In {\jr this} section $\rho_{\Phi,r}$ will stand
for the extended premetric introduced {\jr at the beginning of the  subsection \ref{Sec_prelim2}}.

\begin{theorem}\label{Thm_fixed-point}
Consider $R$ given as in Definition \ref{Def_SpaceE_F},  and
$\Phi\in \mathcal{E}_F^R$, with $F$ as in Proposition
\ref{Thm_Closure2}. Let $r$ satisfy the above assumptions including
(A$_r$), and such that (C$_{r,R}$) is satisfied. Let $\Upsilon:
{\rm Hom}(F) \rightarrow {\rm Hom}(F)$ be an application. Let $\{f_n\}$ be a
sequence in {\rm Hom}($F$). We assume that there exists $h_0 \in${\rm Hom}($F$)
such that the following conditions hold:
\begin{itemize}
  \item[(i)] $\delta:=\rho_{\Phi,r}(\Upsilon(h_0),h_0)
  <{\mkk A^{-1}}$, where $A=\max(a\beta, bm^{-1}\beta+\beta \gamma^{-1})$. 
  \item[(ii)] $\{\Upsilon^n(h_0)\}_{n\in \mathbb{Z}}$ is bounded on
  every compact of $F$.
\end{itemize}

Assume furthermore that there exists a constant $0<C<1,$ such that,
\begin{equation}\label{Eq_contraction}
 \forall \; (f,g) \; \in {\rm Hom}(F) \times {\rm Hom}(F), \;
\rho_{\Phi,r}(\Upsilon(f),\Upsilon(g)) < C\rho_{\Phi,r}(f,g),
\end{equation}
 then
there exists a unique $h\in {\rm Hom}(F)$ such that $\Upsilon(h)=h, $
which is obtained as a limit in $\mathfrak{T}^{-}(\rho_{\Phi,r})$ of
$\{\Upsilon^n(h_0)\}_{n\in \mathbb{N}}.$

Furthermore, if there exists $n_0 \in \mathbb{N}$ such that $
\Upsilon^{n_0}(h_0) \in \mathbb{H}_{\Phi,r},$ then $h\in
\mathbb{H}_{\Phi,r}.$
\end{theorem}

{\mkk \begin{remark}
Note that  $A> 1$ since $\beta>\gamma$ by assumption (G$_3$).
\end{remark}
}

\begin{proof}
From Proposition \ref{Thm_Closure2}, it suffices to show that
$\{\Upsilon^n(h_0)\}_{n\in \mathbb{N}}$ is $\rho_{\Phi,r}^+$-Cauchy
for any $h_0 \in $ Hom($F$) satisfying (i). For simplifying the
notations we denote by $\rho_n^m$ the quantity
$\rho_{\Phi,r}(\Upsilon^m(h_0),\Upsilon^n(h_0))$, for $m\geq n$.
Note that since $\delta $ is finite, by recurrence and using the
contraction property \eqref{Eq_contraction} we can show that all the
quantities $\rho_n^m$ are finite as well.

By using the relaxed inequality \eqref{Eq_relaxed_ineq2}, the
contraction property \eqref{Eq_contraction} and the definition of
$A$ in condition (i), it {\mkk is easy  to obtain} for any integers $n, m \geq
n+1$, and $k\geq 1$,

\begin{equation}\label{Eq_int21}
\begin{split}
\rho_{\Phi,r}(\Upsilon^{m+k}(h_0),\Upsilon^{n}(h_0)) \leq A \delta &
C^{m+k-1} \rho_{\Phi,r}(\Upsilon^{m+k-1}(h_0),\Upsilon^{n}(h_0))
\\
& + A\delta C^{m+k-1} +
\rho_{\Phi,r}(\Upsilon^{m+k-1}(h_0),\Upsilon^{n}(h_0)),
\end{split}
\end{equation}
which leads to,
\begin{equation}\label{Eq_int21}
\rho_n^{m+k} < C^{m+k-1} \rho_n^{m+k-1} +  C^{m+k-1} +
\rho_n^{m+k-1},
\end{equation}
by using $A\delta <1$ from assumption (i) and the notations
specified above.

Let $\epsilon >0$ be fixed. For any $m$ and $n$,  we introduce now
the two-parameters sequence $\{F_k(m,n)\}_{k\in \mathbb{N}}$ defined
by recurrence through,
\begin{equation}\label{Eq_Suite_auxiliaire}
   \left\{ \begin{array}{l}
  F_k(m,n)=C^{m+k-1} F_{k-1}(m,n) + C^{m+k-1}+  F_{k-1}(m,n), \forall \; k\geq 1,\\
  F_0(m,n)=\epsilon.
\end{array} \right.
\end{equation}
When no ambiguity is possible, $F_k$ will simply stand for
$F_k(m,n)$. The role of $m$ and $n$ will be apparent in a moment.

Moreover, from \eqref{Eq_int21}, {\jr it is} easy to show that for any $n$
and $m\geq n+1$,
$$(\rho_n^m \leq F_0(m,n))\Rightarrow (\forall \; k\geq 1, \;  \rho_n^{m+k} \leq F_k(m,n)). $$
As we will see, it suffices to consider $m=n+1$ to prove the
theorem, a choice that we make in what follows. Since $C<1$,
obviously,
$$ \exists \; N_1\; :\; \forall\; n\geq N_1, \rho_n^{n+1}< C^n \delta <\epsilon=F_0(n+1,n),$$
and therefore we get that $\rho_n^{n+1+k} \leq F_k(n+1,n),$ from
what precedes.

The key idea is now to note that if for all $k \geq 0$,
$F_k(n+1,n)\leq 2 \epsilon,$ for $n$ sufficiently big, then the
sequence $\{\Upsilon^n(h_0)\}_{n\in \mathbb{N}}$ is
$\rho_{\Phi,r}^+$-Cauchy for $h_0 \in $ Hom($F$) fulfilling
condition (i). In the sequel we prove that it is indeed the case.

To do so, an easy recurrence shows that,
$$ \forall \; k\in
\mathbb{N}, \; F_k>0, \mbox{and } \{F_k\} \; \mbox{is strictly
increasing}.$$
In particular,
$$ \forall \; k\in
\mathbb{N}, \; F_k \geq \epsilon, $$ and therefore for any $k\geq 1$
and $m$,
\begin{equation}\label{Eq_int22}
\frac{F_k}{F_{k-1}}=C^{m+k-1} +1 + \frac{C^{m+k-1}}{F_{k-1}}\leq
C^{m+k-1}\big(1+\frac{1}{\epsilon}\big) +1.
\end{equation}

Thus by using \eqref{Eq_Suite_auxiliaire} and iterating
\eqref{Eq_int22}, we obtain for any $k\geq 1$,
\begin{equation}\label{Eq_int23}
v_k:=F_k-F_{k-1}< C^{m+k-1} (F_{k-1}+1) \leq C^{m+k-1} \Big\{
\prod_{l=2}^k \Big(C^{m+l-2} (1+\frac{1}{\epsilon})+1\Big)\cdot
\epsilon +1\Big\},
\end{equation}
with the convention $\prod_{l=2}^{k=1} \Big(C^{m+l-2}
\big(1+\frac{1}{\epsilon}\big)+1 \Big)\equiv 1,$ making valid
\eqref{Eq_int23} for $k=1$.

Since $C<1$, then for any $l\in \{2,...,k\}$, and $k\geq 2$,
$C^{m+l-2}\leq C^m$, which leads from \eqref{Eq_int23} to,

\begin{equation}
v_k \leq C^m\Big\{\Big(C^{m+1}\big(1 +\epsilon^{-1} \big)+C
\Big)^{k-1} \cdot \epsilon + C^{k-1} \Big\},
\end{equation}
which is also valid for $k=1$, by simply computing $F_1-F_0$.

Besides,
\begin{equation}
\exists \; N_2 \; :\; \forall \; m\geq N_2, \; a_m:= C^{m+1}
(1+\epsilon^{-1}) +C <1,
\end{equation}
which shows for $m\geq N_2$,

\begin{equation}\label{Eq_int24}
\sum_{j=1}^{j=k} v_j \leq C^m\Big\{\frac{\epsilon}{1-a_m} +
\frac{1}{1-C}\Big\},
\end{equation}
since $C<1$ and $a_m<1$. Now it can be shown,
\begin{equation}
\exists \; N_3\;: \; \forall \; m\geq N_3, \; C^m
\Big\{\frac{\epsilon}{1-a_m} + \frac{1}{1-C}\Big\}\leq \epsilon.
\end{equation}
Fixing $m=n+1,$ we conclude from \eqref{Eq_int24} and the trivial
identity $F_k=\sum_{j=1}^{j=k} v_j +F_0$ that,
\begin{equation}
\forall\; n\in\mathbb{N}, \forall k\in \mathbb{N}, \; (n\geq
\max(N_1,N_2,N_3))\Rightarrow (F_k(n+1,n))\leq 2\epsilon),
\end{equation}
which shows in particular that,
\begin{equation}
\forall\; n\in\mathbb{N}, \forall k\in \mathbb{N}, \; (n\geq
\max(N_1,N_2,N_3))\Rightarrow (\rho_n^{n+k+1}\leq 2\epsilon).
\end{equation}
We have thus proved that $\{\Upsilon^n(h_0)\}_{n\in \mathbb{N}}$ is
$\rho_{\Phi,r}^{+}$-Cauchy for any $h_0 \in$ Hom($F$) fulfilling
conditions (i) and (ii), and thus by Proposition \ref{Thm_Closure2},
$h:=\underset{n\rightarrow\infty}\lim \Upsilon^n(h_0)$ exists in
$\mathfrak{T}^{-}(\rho_{\Phi,r}).$

It can be shown furthermore that, for every $h_0 \in $ Hom($F$),
$\{\Upsilon^n(h_0)\}_{n\in\mathbb{N}}$ is incrementally bounded with
respect to any $n_0 \in \mathbb{N},$ due to the contraction property
and the fact  that $\{\Upsilon^n(h_0)\}_{n\in\mathbb{N}}$ is
$\rho_{\Phi,r}^+$-Cauchy. {\jr Consequently}, if there exists $n_0$ such
that $\Upsilon^{n_0}(h_0) \in \mathbb{H}_{\Phi,r}$ then by applying
Proposition \ref{Thm_Closure2} again (last part) we obtain $h\in
\mathbb{H}_{\Phi,r}$.

\end{proof}

\section{A conjugacy theorem and the generalized spectrum of the Koopman
operator}\label{Sec_main_thm}

\subsection{The conjugacy theorem}\label{Sec_main_thm'}
We prove in this section the main result of this article, {\it i.e}
the conjugacy Theorem \ref{THM_conjug}. To do so, we need further
preliminary tools and notations that we describe hereafter. In {\jr this} 
section we assume the previous assumptions on $r$ (see {\jr subsection
\ref{Sec_prellim1}}) including the condition (A$_r$). As in
Section \ref{Sec_Fixed_Pt_Thm} the premetric $\rho_{\Phi,r}$ will
stand for the extended premetric introduced in subsection
\ref{Sec_prelim2}. In what follows, we endow again Hom($F$) with a
topology $\mathfrak{T}^{-}(\rho_{\Phi,r})$ where $r$ and $R$ satisfy
the condition (C$_{r,R}$) of {\jr subsection
\ref{Sec_prellim1}},
except that here $\Phi$  belonging to $\mathcal{E}_F^R$ is not
arbitrary and has to solve a generalized eigenvalue problem to
handle the conjugacy problem; cf Theorem \ref{THM_conjug}.

For any self-mapping $f$ of $F$, we introduce the following
$r$-Lipschitz constant,
\begin{equation}
\lambda_r(f):=\sup\Big\{\frac{r(\|f(x)-f(y)\|)}{r(\|x-y\|)}, \; x,
y\in F, \; x\neq y\Big\},
\end{equation}
which can {\jr be infinite}.{\mkk This quantity is clearly a direct
extension of  the {\mkk classical notion of} Lipschitz
constant, Lip($f$),  of a function $f$ on normed vector space where the norm, $\|\cdot\|$, has been substituted by the subadditive map $r(\|\cdot\|)$. }

\begin{definition}
We denote by $\mathbb{L}_r(F)$ the set of homeomorphisms of $F$ such
that $\lambda_r(f),$ and $\lambda_r(f^{-1})$ exist. Such an
homeomorphism is called an $r$-{\mk Lipeomorphism} of $F$.
\end{definition}

We will also need the following concept of generalized eigenvalue of
the Koopman operator, outlined in the introduction.

\begin{definition}\label{Def_Koopman_gene_eigen}
Let $f\in$ Hom($F$), and $U_f$ its Koopman operator with domain
$\mathcal{E}_F^R$. A generalized eigenvalue of $U_f$ is any $\lambda
\in \mathbb{R}$ such that,
\begin{equation}\label{Eq_eigen_gene}
 \exists \; \Phi \in \mathcal{E}_F^R \;: \;  U_f(\Phi)\geq \lambda \Phi, 
 \end{equation}
where in case of existence, $\Phi$ is the  corresponding
generalized eigenfunction.
\end{definition}

{\mkk \begin{remark}
Note that if $f \in \mathbb{H}_{\Phi,r}$  then $U_f(\Phi) \in  \mathcal{E}_F^R$ for any generalized eigenfunction $\Phi$, from Proposition \ref{Thm_group} and Remark \ref{Rem_Koopman_stable}. The case 
of equality makes thus sense in \eqref{Eq_eigen_gene}, justifying {\it de facto} the terminology. When $f\in$ Hom($F$)$\backslash \mathbb{H}_{\Phi,r}$, note that the generalized eigenvalue problem \eqref{Eq_eigen_gene}
may still exhibit solutions in some appropriate  space $\mathcal{E}_F^R$ such that $U_f(\Phi) \in  \mathcal{E}_F^R$, as illustrated in subsection \ref{SEC_example}.
\end{remark}}

Lastly, given $f$ and $g$ in Hom($F$), we introduce the following
classical {\it conjugacy operator},
\begin{equation}\label{L_fg}
\mathcal{L}_{f,g} : \left\{
\begin{array}{l}
\mbox{Hom}(F) \to \mbox{Hom}(F)\\
h\to \mathcal{L}_{f,g}(h):=f\circ h\circ g^{-1}.
\end{array} \right.
\end{equation}

We are now in {\jr a} position to state and prove the main theorem of this
article, a conjugacy theorem which is conditioned to a generalized
eigenvalue problem of the related Koopman operators.

\begin{theorem}\label{THM_conjug}
Given $f$ and $g$ in {\rm Hom}($F$) where $F$ is as in Proposition
\ref{Thm_Closure2}, assume that there exist a growth function $R$
 given as in Definition \ref{Def_SpaceE_F} and a function $r$
satisfying the above assumptions including (A$_r$),  such that
(C$_{r,R}$) is satisfied; and such that the following conditions are
fulfilled,
\begin{itemize}
\item[(a)] $f,g\in \mathbb{L}_r(F)$,
\item[(b)] there exists $\alpha >1$ and a common generalized eigenfunction $\Phi
\in \mathcal{E}_F^R$ of the Koopman operators $U_f$ and $U_g$, which
solves the following generalized eigenvalue problem,
\begin{equation}\label{System_cohomo}  \mathbf{P}_{\alpha}: \left \{
\begin{array}{l}
U_f(\Phi)\geq \alpha \lambda_r(f)\Phi, \\
U_g(\Phi) \geq \alpha \lambda_r(g)  \Phi.\\
\end{array} \right.
\end{equation}
\end{itemize}

Under these conditions, for any $\Phi$ solving
\eqref{System_cohomo}, assume further that there exists an
homeomorphism $h_0$ of $F$ satisfying the following properties:
\begin{itemize}
\item[(i)] $\delta:=\rho_{\Phi,r}(\mathcal{L}_{f,g}(h_0),h_0)<{\mkk A^{-1}}$, where $A=\max(a\beta, bm^{-1}\beta+\beta \gamma^{-1})$. 
\item[(ii)] $\{\mathcal{L}_{f,g}^n(h_0)\}_{n\in\mathbb{Z}}$ is bounded on every compact
of $F$.
\end{itemize}

Then $f$ and $g$ are conjugated by a unique element $h$ of {\rm Hom}($F$),
which is the limit in $\mathfrak{T}^{-}(\rho_{\Phi,r})$ of
$\{\mathcal{L}_{f,g}^n(h_0)\}_{n\in\mathbb{N}}$. Furthermore, if there exists $n_0\in \mathbb{N}$ such that
$\mathcal{L}_{f,g}^{n_0}(h_0) \in \mathbb{H}_{\Phi,r}$, then
$h\in\mathbb{H}_{\Phi,r}$.
\end{theorem}

\begin{proof}

Let $f$ and $g$ be two homeomorphisms of $F$. Let the function $r$
and the growth function $R$ be such that the conditions (a) and (b)
of  Theorem \ref{THM_conjug} are satisfied. Let $\Phi \in
\mathcal{E}_F^R$ be a solution of \eqref{System_cohomo}. We endow
then Hom($F$) with the topology $\mathfrak{T}^-(\rho_{\Phi,r})$ for
such  a $\Phi$.

Since the existence of a solution in Hom($F$) to the conjugacy
problem is equivalent to the existence of a fixed point in Hom($F$)
of $\mathcal{L}_{f,g}$, it suffices from Theorem
\ref{Thm_fixed-point} to examine if $\mathcal{L}_{f,g}$ is a
contraction in the premetric $\rho_{\Phi,r}$. To do so, we need to
estimate $|\mathcal{L}_{f,g}(h_1)\circ
(\mathcal{L}_{f,g}(h_2))^{-1}|_{\Phi,r} $  as well as
$|(\mathcal{L}_{f,g}(h_1))^{-1} \circ
\mathcal{L}_{f,g}(h_2)|_{\Phi,r} $ for all $h_1, h_2 \in$ Hom($F$).
Simple computations show that for all $x\in F$,
\begin{equation}
\frac{r(\|f\circ h_1\circ h_2^{-1}\circ f^{-1}(x) -x\|)}{\Phi(x)}
\leq \lambda_r(f) \frac{r(\|h_1 \circ h_2^{-1}\circ
f^{-1}(x)-f^{-1}(x)\|)}{\Phi(x)},
\end{equation}
which leads to,
\begin{equation}
|\mathcal{L}_{f,g}(h_1)\circ (\mathcal{L}_{f,g}(h_2))^{-1}|_{\Phi,r}
\leq \lambda_r(f)\cdot \underset{u\in
F}\sup\Big(\frac{\Phi(u)}{\Phi(f(u))}\Big)\cdot |h_1\circ
h_2^{-1}|_{\Phi,r}\;.
\end{equation}
Similar computations show,
\begin{equation}
|(\mathcal{L}_{f,g}(h_1))^{-1}\circ \mathcal{L}_{f,g}(h_2)|_{\Phi,r}
\leq \lambda_r(g)\cdot \underset{u\in
F}\sup\Big(\frac{\Phi(u)}{\Phi(g(u))}\Big)\cdot |h_1^{-1}\circ
h_2|_{\Phi,r}\;,
\end{equation}
and since $\Phi$ solves the generalized eigenvalue problem
$\mathbf{P}_{\alpha}$, we get for all $u\in F$,
\begin{equation}
\lambda_r(f) \frac{\Phi(u)}{\Phi(f(u))} \leq \frac{1}{\alpha}<1,
\mbox{ and } \lambda_r(g) \frac{\Phi(u)}{\Phi(g(u))} \leq
\frac{1}{\alpha}<1,
\end{equation}
which allows us to conclude that,
\begin{equation}
\rho_{\Phi,r}(\mathcal{L}_{f,g}(h_1), \mathcal{L}_{f,g}(h_2)) \leq
\frac{1}{\alpha} \cdot \rho_{\Phi,r}(h_1,h_2),
\end{equation}
for all $h_1$  and $h_2$ in Hom($F$), {\it i.e.} the conjugacy
operator $\mathcal{L}_{f,g}$ is a contraction for the premetric
$\rho_{\Phi,r}$. The rest of the assumptions (i) and (ii) of Theorem
\ref{THM_conjug} are just a translation of the ones used in Theorem
\ref{Thm_fixed-point}, and thus by using this last theorem the proof
of the present one is easily achieved.
\end{proof}

\begin{remark}\label{Rem_control}
{\jr This theorem provides conditions for  the conjugacy $h$,  when it exists, to lie in $\mathbb{H}_{\Phi,r}.$}  This means in such a case that $h$ satisfies
some behavior at infinity prescribed by $\Phi$, which has in turn to
solve $\mathbf{P}_{\alpha},$ a spectral problem related to the
Koopman operators $U_f$ and $U_g$,  which involves structural
constants of $f$ and $g$: $\lambda_r(f)$ and $\lambda_r(g)$ as
introduced above. This aspect could be of interest in control
theory.
\end{remark}

\begin{remark}\label{Rem_bounded_every_compact}
Condition of type $(i)$ in Theorem \ref{THM_conjug} is often met in
a stronger form when dealing with (local) conjugacy problems that
arise around a fixed point. Indeed it is often required that the
convergence of the sequence $\{f^{n} g^{-n}\}$ holds in C$^0$
provided that $f$ and $g$ are tangent to sufficiently high order; cf.  {\it e.g.} \cite[Lemma 3 p. 95]{TOPOChaperon86}.
\end{remark}

\subsection{An illustrative example}\label{SEC_example}

To simplify we set $E=\mathbb{R}$ and we consider $F=[0,+\infty)$
which fulfills the conditions of  Theorem \ref{THM_conjug}. We
consider furthermore $r(x)=x$ and $\Phi(x)=R(x)=\sqrt{ x}+1$ (for
$x\in F$) that are subadditive function on $F$. Note that such an
$R$ satisfies the {\mkk conditions of  subsection \ref{Sec_prellim1},
and that $\Phi \in \mathcal{E}_F^R$ with $m=1$ in (G$_1$)  and,  for instance,  $\beta =2$ and $\gamma=1/2$ for (G$_3$).} 

 Note also that $r$ as continuous function satisfies (A$_r$) (and (S)) and
fulfill all the other standing assumptions. The cross condition
(C$_{r,R}$) of {\mkk  subsection \ref{Sec_prellim1} is trivially fulfilled, since for instance:}
$$ 1+\sqrt{u}\leq u+\frac{5}{4}, \; \forall \; u\in \mathbb{R}^+.$$

We consider the following dynamics on $F$: ${\mkk f(x)}=\eta x$ with
$0<\eta<1$ and $g(x)=\eta x+\varphi(x)$, where $\varphi\in
C^1(\mathbb{R}^+, \mathbb{R}^+)$ with compact support is such that
$g$ is an homeomorphism of $F$; $\varphi$ will be further
characterized in a moment. Other conditions on $\varphi$ will be
imposed below. Lastly note that $\lambda_r(f)=\mbox{Lip}(f)$ for
$r(x)=x$. Since $\eta<1$, there exists $\epsilon>0$ such that,
\begin{equation}\label{Eq_cohom_Example1}
\frac{\sqrt{\eta}}{\eta}> 1+\epsilon.
\end{equation}
Now since $\eta(1+\epsilon) \leq \sqrt{\eta}<1$, we get for all
$x\in F,$
\begin{equation}
\alpha \mbox{Lip}(f)\Phi(x)=\eta(1+\epsilon)(\sqrt{x}+1)\leq \sqrt{\eta x}+1=\Phi(f(x)),
\end{equation}
with $\alpha=1+\epsilon,$ which shows that $\Phi$ is a
 generalized eigenfunction in $\mathcal{E}_F^R$ of $U_f$ with eigenvalue $\lambda= (1+\epsilon) \mbox{Lip}(f)$ in {\jr this} particular context.

From (\ref{Eq_cohom_Example1})  we get,

\begin{equation}\label{Eq_int26}
\exists \; \epsilon_2>0, \;: \;
\frac{\sqrt{\eta}}{\eta+\epsilon_2}>(1+\epsilon), \mbox{ and }
\epsilon_2 < \eta.
\end{equation}

Besides, for such an $\epsilon_2$, there exists a function
$\varphi\in C^1(\mathbb{R}^+, \mathbb{R}^+)$ with compact support
such that,
\begin{equation}\label{Eq_control_derivé_phi}
 \mbox{Lip}(g)=\underset{x\in F}\max|\eta+\varphi'(x)|=\eta +\epsilon_2,
\end{equation}
and such that  $g$ is still an homeomorphism of $F$ (since
$\epsilon_2 <\eta $).

From \eqref{Eq_control_derivé_phi}  and \eqref{Eq_int26} we get now,
\begin{equation}
\begin{split}
(1+\epsilon)\mbox{Lip}(g)\Phi(x)&=(1+\epsilon)(\eta+\epsilon_2)(\sqrt{x}+1)  \\
& \leq \sqrt{\eta x}+1 \leq \sqrt{\eta x+\varphi(x)}+1=\Phi(g(x)),
\end{split}
\end{equation}
which shows that $\Phi$ is a generalized eigenfunction in
$\mathcal{E}_F^R$ of $U_g$ with $\lambda= (1+\epsilon)
\mbox{Lip}(g)$. We are then left with a common eigenfunction of
$U_f$ and $U_g$ satisfying $\mathbf{P}_{\alpha}$ with
$\alpha=1+\epsilon$. From our assumptions, {\mkk it is easy} to check
furthermore that $f$ and $ g$ belong to $\mathbb{L}_r(F)$, and
therefore conditions (a) and (b) of Theorem \ref{THM_conjug} are
satisfied in {\jr this}  particular setting. Recall from the proof of
Theorem \ref{THM_conjug} that these conditions ensure the
contraction of the conjugacy operator $\mathcal{L}_{f,g}$ for the
premetric $\rho_{\Phi,r}$.

To apply Theorem \ref{THM_conjug} we have now to check the remaining
conditions (i) and (ii) {\mkk for a convenient $h_0\in{\rm Hom}(F)$}.  Let us take $h_0=g$.  We check first condition (i). {\mkk In that respect, we} have  to
estimate $|\mathcal{L}_{f,g}(h_0)\circ
h_0^{-1}|_{\Phi,r}=|fg^{-1}|_{\Phi,r}$ and
$|(\mathcal{L}_{f,g}(h_0))^{-1}\circ
h_0|_{\Phi,r}=|f^{-1}g|_{\Phi,r}$ since $h_0=g$.
Note that there exists {\jr $\psi \in C^1(\mathbb{R}^+,\mathbb{R}^+)$}
such that $g^{-1}(x)=x/\eta + \psi(x),$ and $\psi $ {\jr has} compact
support. 
 Note {\mkk also} that we can find $\varphi$ and thus $\psi$ such that,
$$\nu:=\max\Big(\underset{x\in F}\max|\varphi|,\underset{x\in F}\max|\psi|\Big)<\eta A^{-1},$$
without violating (\ref{Eq_control_derivé_phi}) and thus having
$\Phi$ still satisfying $\mathbf{P}_{\alpha}$. For such a choice,
$\eta\nu <\eta^2A^{-1}<A^{-1}$ since $\eta<1$, and in particular,
\begin{equation}
|fg^{-1}|_{\Phi,r} =\underset{x\in F}\sup \Big(\frac{|\eta
\psi(x)|}{\sqrt{x}+1}\Big)<A^{-1},
\end{equation}
and,
\begin{equation}
|f^{-1}g|_{\Phi,r} =\underset{x\in F}\sup
\Big(\frac{|\varphi(x)|}{\eta(\sqrt{x}+1)}\Big)<A^{-1},
\end{equation}
which allows us to conclude that condition (i) is checked with
$h_0=g$.

 Finally let us check condition (ii). Since $\varphi$ and $\psi$ {\jr have compact supports}  then it can be
shown that for $h_0=g$, the sequence
$\{\mathcal{L}_{f,g}^n(h_0)\}_{n\in\mathbb{Z}}$ is bounded on every
compact {\mk subset} of $F$. We leave the details to the reader. We can thus
apply Theorem \ref{THM_conjug} to conclude that $f$ and $g$ are
conjugated which was of course obvious for Lip($\varphi$)
sufficiently small, from a trivial application of the global
Hartman-Grobman theorem, cf. for instance \cite{TOPOIrwin}.

This modest example is just intended to illustrate some mechanisms
of the approach developed in this article. Of course, further
investigations are needed with respect to the existence of solutions
of the generalized eigenvalue problem $\mathbf{P}_{\alpha}$ in
spaces of type $\mathcal{E}_F^R$, for more general homeomorphisms.
We postpone this difficult task for a future work, discussing in the
next subsection some related issues.

\subsection{Generalized {\mkk eigenfunctions} of the Koopman operator}\label{SEC_COHOM_discussion}
We describe here two possible approaches to examine the generalized
eigenvalue problem $\mathbf{P}_{\alpha}$, the first one is based on
Schr$\mbox{\"{o}}$der equations and the second one is based on
cohomological equations. The point of view retained is based on
functional equations techniques coming from different part of that
literature where we emphasize the overlapping. In both cases, by
shortly reviewing the existing results, we provide hereafter
conditions under which the generalized eigenvalue problem
$\mathbf{P}_{\alpha}$ may possess continuous solutions, without
being able to specify --- in a general setting --- if $\Phi$ can
live in {\mkk some space} $\mathcal{E}_F^R$. {\jr We elaborate on this point here with the intent of gathering some results related to our problem which are found dispersed in the literature.} 

Note that in the sequel, we will focus more precisely on
the generalized eigenvalue problem for $U_f$ and not
$\mathbf{P}_{\alpha}$ itself, in order to exhibit already the main
issues for the associated single existence problem of a generalized
eigenfunction.

\subsubsection{Approach based on Schr$\mbox{\"{o}}$der
equations}\label{Sect_Schroder} We recall first some background
concerning Schr$\mbox{\"{o}}$der equation. Here $E$ denotes a real
or complex normed vector space, $\mathcal{H}$ denotes an
{\mkk arbitrary space  of $\mathbb{R}-$ or $\mathbb{C}-$valued functions on $E$, and $\mathcal{F}$ denotes some space of self mappings of $E$}. Let {\mkk $f$ be in $\mathcal{F}$}, then
the Schr\"{o}der's equation in $\mathcal{H}$ is the
equation of unknown $\Psi$ ({\mkk to be found in $\mathcal{H}$}):
\begin{equation}
\Psi \circ f =\lambda \cdot \Psi,
\end{equation}
for some $\lambda$. It is the equation related to the spectrum of
the Koopman operator\footnote{The Koopman operator is also known as
the composition operator in other fields \cite{TOPOCowen_book}.}
U$_f: \Psi \mapsto \Psi\circ f$, in the space $\mathcal{H}$. The
properties of this spectrum are closely related to the function $f$
as well the space $\mathcal{H}$. The Schr$\mbox{\"{o}}$der's
equation has a long history and has been extended and studied in
various settings. In the early 1870s, Ernst Schr$\mbox{\"{o}}$der
\cite{TOPOSchroder} studied this type of functional equation in the
complex plane for the composition operator, for $\Psi(z) = z^2$ and
$f(z) = z + 1$. The functional equation named after him is $\Psi
\circ f = \lambda \Psi$ where $f$ is a given complex function, and
the problem consists of finding $\Psi$ and $\lambda$ to satisfy the
equation, {\it i.e.} an eigenvalue problem for $U_f$. An important
part of the results in the literature are devoted to contexts where
$f$ is a function mapping the unit disc in the complex plane onto
itself initiated by the seminal work of Gabriel Koenigs in 1884
\cite{TOPOkoenigs}. The reader may consult \cite{TOPOClahane} or \cite{TOPOShapiro}  
 with references therein, for a recent account
about this part of the literature. This functional problem has also
been considered historically for {\jr maps}  of the half-line
\cite{TOPOKuczma0} or more general Banach spaces in the past decades
\cite{TOPOWalorski08}. We mention lastly, that the
Schr$\mbox{\"{o}}$der equation is sometimes encountered under the
form of the Poincar\'e functional equation \cite{TOPOJulia, TOPOKuczma0} and arises in various applications such {\jr as} iterated
function theory \cite{TOPODenjoy, TOPOKuczma}, branching process
\cite{TOPOSeneta} or dynamical systems {\mkk theory \cite{TOPObelitskii_Schroder,TOPOYoccoz}}.

Naturally, the generalized eigenvalue problem $\mathbf{P}_{\alpha}$
can be related to Schr$\mbox{\"{o}}$der equations. Indeed, if the
following  Schr$\mbox{\"{o}}$der equation,
\begin{equation}\label{Schroder_eq}
\Psi \circ f = \alpha \lambda_r(f) \cdot \Psi,
\end{equation}
has a solution $\Psi:F \rightarrow F$ for $\alpha >1$, then the
generalized eigenvalue problem  for $U_f$ has an obvious solution,
provided that $\Phi(\cdot)=\|\Psi(\cdot)\| \in \mathcal{E}_F^R$ and
$f\in \mathbb{L}_r(F)$. It is known that such an equation can be
solved for particular domain $F$ and particular space of {\mkk functions}
over $F$ such as Hardy spaces; see \cite{TOPOCowen_book}. It is
interesting to note that most of the results typically {\mkk require} some
{\jr compactness}  assumptions of the Koopman operator $U_f$ which involve
that $f$ possesses at least a fixed point in  $F$; {\mkk cf.
 \cite{TOPOCaughran_Schwartz,TOPOShapiro_al}}, cf. also  
\cite[Theorem 5.1]{TOPOClahane} for extensions of results of
\cite{TOPOCaughran_Schwartz}.

There exist other techniques coming from functional analysis rather
than complex analysis, to deal with $\mathbf{P}_{\alpha}$ from the
point of view of Schr$\mbox{\"{o}}$der equations. It consists of
considering a more general type of Schr$\mbox{\"{o}}$der equations
where the unknown is a map $\Psi:F \rightarrow E$ aiming to satisfy,
\begin{equation}\label{Schroder_eq2}
U_f(\Psi) =  A \circ  \Psi,
\end{equation}
where $A$ is a linear map of $E$. If we assume that $A$ is
invertible, and that there exist a C$^0$-functional
$\mathcal{N}:E\rightarrow \mathbb{R}$ with a constant
 $\mathfrak{m}>0$ satisfying,
\begin{equation}\label{schroder_functional_Walorski}
\forall\; x\in F, \quad \mathcal{N}(f(x))\geq \mathcal{N}(x)
+\mathfrak{m},
\end{equation}
then  \cite[Theorem 2.1]{TOPOWalorski08} {\jr permits to conclude  the}
existence of a {\it continuous nonzero} solution of
(\ref{Schroder_eq2}). Note that if $F\subset E\backslash \{0\}$ and
$\|f(x)\| \leq \kappa \|x\|$ on $F$, with $\kappa \in (0,1)$ then
the functional inequality \eqref{schroder_functional_Walorski} is
satisfied on $F$ by simply taking $\mathcal{N}(x):=-\log\|x\|$ and
$\mathfrak{m}:=-\log(\kappa)$.

Now by assuming $A$ invertible,
$$ \forall \; \xi\in E,\; \|\xi\| \leq \|A^{-1}\| \|A \xi\|,$$
$$\mbox{{\it i.e.}}, \;\|A^{-1}\|^{-1}\leq \underset{\xi\in E-\{0\}}\inf\frac{\|A\xi\|}{\|\xi\|},$$
and thus if,
\begin{equation}\label{A-1spectrum_ineq}
\|A^{-1}\|^{-1}\geq \alpha \lambda_r(f),
\end{equation}
and there exists a couple $(\mathfrak{m},\mathcal{N})$ satisfying
\eqref{schroder_functional_Walorski}, we obtain that
$\Phi(\cdot):=\|\Psi(\cdot)\|$ with $\Psi$ a solution of
(\ref{Schroder_eq2}), is a solution of $U_f(\Phi) \geq \alpha
\lambda_r(f) \cdot \Phi$. Thus in order to have a solution of the
generalized eigenvalue problem  for $U_f$, we only need to know
about the asymptotic behavior of $\Phi$. However, the examination of the growth of an eigenfunction $\Psi$
solution of (\ref{Schroder_eq2}) is a difficult task in general, \footnote{
 {\it e.g.} \cite{TOPOBour_Shapiro, TOPOShapiro_al} or \cite{TOPOBanas} for
particular cases   related to the standard
Schr\"oder equation (\ref{Schroder_eq}).} which renders the
generalized eigenvalue problem  for $U_f$ and thus the functional
problem $\mathbf{P}_{\alpha}$ introduced here non-trivial to solve
in general.

To conclude, we emphasize that the Schr$\mbox{\"{o}}$der equation is
related to Abel's functional equation \cite{TOPOAbel}, which is a
well known functional equation often presented into the form
$\varphi(f(x))=\varphi(x)+1$, where $\varphi :X \rightarrow
\mathbb{C}$ is an unknown function and $f:X \rightarrow X$ is a
given continuous mapping of a topological space $X$
\cite{TOPOKuczma0}.\footnote{The link between the both is trivial
when $X\equiv\mathbb{C}$ where every solution $\varphi$ of the Abel
equation leads to a solution $\Psi:x\mapsto\exp(\log(\lambda)
\varphi(x))$ for every $\lambda
>0$ of the Schr$\mbox{\"{o}}$der equation. In such a case the spectrum of $U_f$ contains $(0,+
\infty)$.} {\jr If the Abel equation
possesses a continuous solution, then  it provides a continuous
solution to (\ref{schroder_functional_Walorski}). Thus  Abel's equation has a  central role in our approach. As the next subsection indicates, there  are  deep connections between both equations and  our functional
problem.}

\subsubsection{Approach based on cohomology equations}\label{Sect_Cohomology}
Even if, to the best of the knowledge of the authors, Theorem
\ref{THM_conjug} exhibits new relations between the existence of
a conjugacy between two homeomorphisms and the spectrum of the
related Koopman operators, relations between conjugacy problems
and functional equations are far to be new.  They arise classically
under the form of the {\it Livshitz cohomology equation}
\cite{TOPOLivschitz, TOPOLivschitz2}, $\phi=\Phi\circ f-\Phi$, where
$f:\M\rightarrow \M$ is a dynamical system of some manifold $\mathfrak{M}$;
$\phi:\M \rightarrow \mathbb{R}$, a given function,  and $\Phi$ maps
$\M$ into $\mathbb{R}$ or a multidimensional space; see for instance
\cite{TOPOBanyanga_Llave, TOPOBanyanga_Llave2,TOPOLMM86, TOPOkat}.
As pointed by Livshitz \cite{TOPOLivschitz, TOPOLivschitz2} the
existence of a continuous solution $\Phi$ strictly depends on the
dynamics generated by $f$ and the topological as well as geometrical
properties of $\M$. For instance if we consider the particular case
of the Abel equation, $\Phi(f(x))=\Phi(x)+1$, there is no continuous
solution if $\M$ is compact, since if such solution would exist,
$\Phi(f^n(x))=\Phi(x)+n$, which would be impossible.

In the case of non-compact topological manifold, Belitskii and
Lyubich {\mkk in \cite{TOPObelitskii_Abel1}} have proved the following theorem, that we present in a
slightly less general setting {\mkk than \cite[ Corollary 1.6]{TOPObelitskii_Abel1}}, adapting their statements with
respect to our purpose:

\begin{theorem}(From \cite{TOPObelitskii_Abel1})\label{THM_Belit}
Assume that $\M$ is locally compact and countable at infinity. If
$f:\M\rightarrow \M$ is continuous and injective then the following
statements are equivalent,
\begin{itemize}
\item[(a)] There exists a continuous solution $\varphi: \M \rightarrow
\mathbb{C}$ of the Abel equation, $\varphi(f(x))=\varphi(x)+1.$
\item[(b)] For every continuous functions $p: \M\rightarrow
\mathbb{C}\backslash\{0\}$ and $\gamma: \M\rightarrow \mathbb{C}$
there exits a continuous solution $\varphi :\M \rightarrow \mathbb{C}
$ of
\begin{equation}\label{Belitskii-eq}
\varphi(f(x))=p(x)\varphi(x)+\gamma(x).
\end{equation}
\item[(c)] Every compact subset of $\M$ is wandering for $f$.
\end{itemize}
\end{theorem}

In the above theorem, a compact set $K \subset \M$, is qualified to
be wandering if there exists an integer $\nu \geq 1$ such that
$$ f^n(K)\cap f^m(K)=\emptyset \quad (n-m \geq \nu),$$
in particular such a dynamical system $f$ is fixed-point free and
periodic-point free, which  {\mkk is for instance} consistent with {\mkk dynamical
restrictions} imposed by any solution to the functional problem
$\mathbf{P}_{\alpha}$ in the case $\lambda_r(f)\geq
1.$\footnote{Indeed, if there exists a generalized eigenfunction
$\Phi$ of $U_f$ and a periodic orbit of $f$ of period $p$ emanating
from some $x^*$, then by repeating $p$-times the change of variable
$x\leftarrow f(x)$ in $\Phi (f(x)) \geq \alpha \lambda_r(f)
\Phi(x)$, we deduce necessarily that $(\alpha \lambda_r(f))^p \leq
1$ (since $\Phi>0$), which imposes that $f$ cannot possess such a
periodic orbit in the case $\lambda_r(f)\geq1$.}

This  theorem provides however an incomplete {\mkk answer} to the
problem $\mathbf{P}_{\alpha}$. For instance, Theorem \ref{THM_Belit}
shows that for $f$ satisfying condition $(c)$ above, there exits a
solution of the equation $\varphi(f(x))=\lambda \varphi(x)$ with
$|\lambda| \geq \alpha \lambda_r(f)$, and therefore a solution of
the generalized eigenvalue problem for $U_f$; obtained by taking its
module, $\Phi(\cdot):=|\varphi(\cdot)|$. The missing step for having
a full solution of that problem is still the knowledge of the
asymptotic behavior of $\Phi$, which is {\mkk also} a difficult {\mkk property to derive for
solutions of cohomology equations}; {\it e.g.} \cite{TOPOBanyanga_Llave2}.

\begin{remark}\label{Importance_F}
It is worth mentioning that in the framework developed in this
article, the subset $F$ plays a central role in the existence of
solutions of functional problem of type $\mathbf{P}_{\alpha}$. The
reader may consult for instance \cite[\S
5.2]{TOPObelitskii_Schroder},  where it is shown that the Abel equation associated with a
contracting mapping of the  cut plane
$\mathbb{R}^2\backslash(-\infty, 0]$, has a real-analytic solution,
whereas the same map considered on the whole plane leads to an Abel
equation without any continuous solution, since such a map possesses
obviously a fixed point which is excluded by Theorem
\ref{THM_Belit}.
\end{remark}

 \section{Concluding remarks}
 {\mk From the previous section}, whatever $F$ and the approach retained, we can conclude that the
main issue concerns therefore the asymptotic behavior of a possible
eigenfunction of the generalized eigenvalue problem
$\mathbf{P}_{\alpha}$. The regularity of the common eigenfunction is
also an important aspect of the problem, in that respect, the
continuity assumption on $R$ and thus on $\Phi$ could be relaxed
using Remark \ref{Rem_mesurability}. Dynamical properties such as
condition (c) of Theorem \ref{THM_Belit} might play also a  role in
the existence of such eigenfunctions. Note also, that since the
closure of $\mathcal{E}_F^R$ in the compact-open topology
\cite{TOPOhirch} is a closed cone with non empty interior in
$C^0(F,\mathbb{R})$ (cf. Remark \ref{Rem_mesurability}-(a)), it would
be interesting to study the spectral properties of $U_f$ and $U_g$
within an approach of type Krein-Rutman theorem \cite{TOPOSchaefer}.
However, the task is more difficult {\mk in the present context} than usually since $F$ is
{\mk assumed to be unbounded}, $\sigma$-compact and locally compact, {\mk which implies that}
$C^0(F,\mathbb{R})$ has a Fr\'echet structure \cite{TOPOAli} (and not a
Banach one), and {\mk in particular makes} non straightforward an extension of the
classical Krein-Rutman theorem in that context {\mk in order} to analyze the
existence of a principal eigenfunction in $\mathcal{E}_F^R$.

Finally, we have intentionally not considered important dynamical
properties such as {\mkk uniform} hyperbolicity {\mk or its violation} \cite{TOPOkat} that could lead to
other spectral problem than $\mathbf{P}_{\alpha}$; {\mk this last one being presented here at a level of generality which lays the foundations for such enterprise.}
{\mk In that perspective, the cone condition (G$_3$) might be relaxed as suggested in Remark \ref{Rem_mesurability} (b), in order to bound the behavior of the  generalized eigenfunction function $\Phi(x)$ by subadditive functions 
only for large values of $x$ for instance, making thus the corresponding problem {\bf P$_{\alpha}$} more flexible.}

{\mk In summary, the} purpose of the present work was to {\mk introduce}, in a general context, a
framework that makes apparent certain relations between the
spectral theory of dynamical systems and the {\mkk topological} problem of
conjugacy. Likewise, the idea of using some observable --- here a
common eigenfunction of the Koopman operators
--- to build a specific topology to deal with the conjugacy problem
does not seem to be limited to the case of unbounded phase-space.

\section*{Acknowledgments}
The authors acknowledge the anonymous referee
for his or her numerous relevant comments which helped to improve the presentation of this article. This work was partially supported by the National Science Foundation
 grant DMS-1049253 and Office of Naval Research grant
N00014-12-1-0911(MDC).

\medskip
\medskip

\end{document}